\documentclass[psamsfonts,reqno]{amsart}
\usepackage{amssymb}
\usepackage{amscd}

\theoremstyle{plain}
\newtheorem{theorem}{Theorem}[section]
\newtheorem{lemma}[theorem]{Lemma}
\newtheorem{proposition}[theorem]{Proposition}
\newtheorem{corollary}[theorem]{Corollary}

\newtheorem*{BT}{Bergman's Theorem}
\newtheorem*{ET}{Elliott's Lemma}
\newtheorem*{ST}{Schmidt's Theorem}
\newtheorem*{PL}{Pudl\'ak's Lemma}
\newtheorem*{BC}{Corollary}

\theoremstyle{definition}

\newtheorem*{notation}{Notation}
\newtheorem{example}[theorem]{Example}
\newtheorem{problem}{Problem}
\newtheorem{remark}[theorem]{Remark}

\newtheorem*{stat}{\name}
\newcommand{\name}{testing}

\newenvironment{all}[1]{\renewcommand{\name}{#1}\begin{stat}}
                        {\end{stat}}

\newcommand{\qedc}{{\qed}~{\rm Claim~{\theclaim}.}}

\DeclareMathOperator{\Id}{Id}
\DeclareMathOperator{\Idbar}{\overline{Id}}
\DeclareMathOperator{\Idc}{Id_c}
\DeclareMathOperator{\Idcbar}{\overline{Id}_c}
\DeclareMathOperator{\Con}{Con}
\DeclareMathOperator{\Conc}{Con_c}

\DeclareMathOperator*{\HH}{H}

\newcommand{\dnw}{\mathbin{\downarrow}}
\newcommand{\rng}{\mathrm{rng}}
\newcommand{\supp}{\mathrm{supp}}
\newcommand{\cm}{commutative monoid}
\newcommand{\crm}{conical refinement monoid}
\newcommand{\poag}{partially ordered abelian group}
\newcommand{\tw}{\mathbf{2}}
\newcommand{\CC}{\mathbb{C}}
\newcommand{\NN}{\mathbb{N}}
\newcommand{\ZZ}{\mathbb{Z}}
\newcommand{\QQ}{\mathbb{Q}}

\newcommand{\AAA}{\mathcal{A}}
\newcommand{\OO}{\mathcal{O}_2}
\newcommand{\Vtil}{\widetilde{V}}

\newcommand{\LL}{\mathcal{L}}
\newcommand{\VV}{\mathbf{V}}
\newcommand{\FF}{\mathbf{F}}
\newcommand{\vtil}{\widetilde{V}}
\newcommand{\vx}{\mathsf{x}}
\newcommand{\diag}{\mathrm{diag}}

\newcommand{\BB}{\mathbf{B}}
\newcommand{\EE}[1]{\QQ\langle#1\rangle}
\newcommand{\EEp}[1]{\QQ^+\langle#1\rangle}
\newcommand{\GG}[1]{\QQ(#1)}
\newcommand{\GGp}[1]{\QQ^+(#1)}

\def\vv<#1>{\langle#1\rangle}

\begin{document}

\title[Representations of distributive semilattices]
{Representations of distributive semilattices in ideal
lattices of various algebraic structures}

\author{K.R. Goodearl}
\address{Department of Mathematics\\
University of California\\
Santa Barbara, CA 93106\\
U.S.A.}
\email{goodearl@math.ucsb.edu}
\urladdr{http://www.math.ucsb.edu/\~{}goodearl/}

\author{F. Wehrung}
\address{CNRS, ESA 6081\\
D\'epartement de Math\'ematiques\\
Universit\'e de Caen\\
14032 Caen Cedex\\
France}
\email{wehrung@math.unicaen.fr}
\urladdr{http://www.math.unicaen.fr/\~{}wehrung}

\keywords{Distributive semilattice, von Neumann regular
ring, dimension group, complemented modular lattice,
C*-algebra, approximately finite dimensional, direct limit,
compact congruence, maximal semilattice quotient, real rank
zero.}
\subjclass{Primary 06A12, 06C20, 06F20, 16E20, 16E50, 19A49, 19K14,
46L05.}

\begin{abstract}
We study the relationships among existing
results about representations of distributive semilattices
by ideals in dimension groups, von Neumann regular rings,
C*-algebras, and complemented modular lattices. We prove
additional representation results which exhibit further
connections with the scattered literature on these different
topics.
\end{abstract}

\maketitle

\section*{Introduction}

Many algebraic theories afford a notion of \emph{ideal}, and
the collection of all ideals of a given object typically
forms a complete lattice with respect to inclusion. It is
natural to ask which lattices can be represented as a lattice
of ideals for a given type of object. Often, the lattice of
ideals of an object is algebraic, in which case this lattice
is isomorphic to the lattice of ideals of the
(join-) subsemilattice of compact elements. For instance, this
holds for lattices of ideals of rings, monoids, and partially
ordered abelian groups. Hence, lattice representation problems
often reduce to corresponding representation problems for
(join-) semilattices.
For example, to prove that a given
algebraic lattice $L$ occurs as the lattice of ideals of a ring
of some type, it suffices to show that the semilattice of
compact elements of $L$ occurs as the semilattice of finitely
generated ideals of a suitable ring.

We shall be concerned here with representation problems for
\emph{distributive} algebraic lattices, which correspond to
representation problems for distributive semilattices. The
contexts we discuss include congruence lattices,
complemented modular lattices, (von Neumann) regular rings,
dimension groups, and approximately finite dimensional
C*-algebras. All these contexts are interconnected, and a
main goal of our paper is to develop these interconnections
sufficiently to allow representation theorems for
distributive semilattices in one context to be transferred to
other contexts.

Since readers familiar with one of our contexts may not be
fully at home in others, we try to provide full details and
all relevant definitions in the appropriate sections of the
paper. While the reader may encounter some undefined
concepts in this introduction, we hope that the flavor of
the results discussed will come through nonetheless on a
first reading. All the required concepts will be made
precise later in the paper.

Typical representation results for distributive semilattices
include the following:

\begin{ST}
Every finite distributive lattice is isomorphic to
the semilattice of compact congruences of some complemented
modular lattice.
\qed
\end{ST}

This is a result of E.T. Schmidt \cite{Schm84}. It is
probably the earliest representation result of distributive
semilattices by complemented modular lattices. Further
lattice-theoretical representation results are discussed in
\cite{GrSc}, mainly in relation
with the \emph{Congruence Lattice Problem}, that asks whether
every distributive algebraic lattice is isomorphic to the
congruence lattice of a lattice.

A stronger version of Schmidt's Theorem follows from a
result of G.M. Bergman~\cite{Berg86}:

\begin{BT}
Let $L$ be a distributive algebraic lattice with
only countably many compact elements, and $K$ any field.
Then there exists a locally matricial
$K$-algebra $R$ of countable dimension whose lattice of
two-sided ideals is isomorphic to~$L$. If, in addition, the
greatest element of $L$ is compact, then one can choose $R$
unital.\qed
\end{BT}

According to the abovementioned correspondence between
semilattices and algebraic lattices, this can also be
formulated as follows: \emph{Every countable distributive
$0$-semilattice $S$ is isomorphic to the semilattice of
finitely generated two-sided ideals in some locally matricial
algebra $R$ of countable dimension. If, in addition, $S$ has
a largest element, then one can choose $R$ unital.}

Locally matricial algebras are, in particular, regular
rings, and the finitely generated right ideals of any
regular ring $R$ form a sectionally complemented modular
lattice, $\LL(R)$. Further, the semilattice of finitely
generated two-sided ideals of $R$ turns out to be isomorphic
to the semilattice of compact congruences of $\LL(R)$, see
Proposition~\ref{P:ManyIso} (\emph{cf.}
\cite[Corollary 4.4]{Wehr2}). Hence, Bergman's Theorem yields
the following result:

\begin{BC}
Any countable distributive $0$-semilattice is
isomorphic to the semilattice of compact congruences of some
sectionally complemented modular lattice.\qed
\end{BC}

The $\aleph_1$ version of Bergman's Theorem is still open
(see the discussion around Problem~\ref{Pb:aleph1Ring} in
Section~\ref{S:Probs}).
The second author has shown that the
$\aleph_2$ version has a
\emph{negative} answer (see \cite{Wehr2}). A precursor to Bergman's
Theorem was obtained by K.H. Kim and F.W. Roush, who proved
that any finite distributive lattice is isomorphic to the
lattice of two-sided ideals of some unital locally matricial
algebra of countable dimension, see \cite[Corollary to
Theorem 4]{KR}. In view of the connections discussed above,
this result is already sufficient to yield Schmidt's Theorem.

An interesting representation result for distributive
semilattices of arbitrary size was proved by P. Pudl\'ak, see
\cite[Fact 4, p. 100]{Pu}:

\begin{PL}
Every distributive semilattice is the direct
union of all its finite distributive subsemilattices.
\qed
\end{PL}

These results are similar in spirit to representation
results in other fields of mathematics, that were proved
completely independently.

We start with G.A. Elliott, who classified countable direct
limits of locally matricial algebras by an invariant
equivalent to their ordered
$K_0$ groups, see \cite[Theorem~4.3]{Elli76}
(\emph{cf.} \cite[Theorem 15.26]{Gvnrr}). Elliott's initial
result towards the question of which ordered groups appear
in this classification \cite[Theorem 5.5]{Elli76} can be
phrased as follows:

\begin{ET}
Let $G$ be the direct limit of a countable
sequence of simplicial groups, and let $K$ be a field. Then
there exists a locally matricial $K$-algebra $R$ of
countable dimension such that $K_0(R)\cong G$. If, in
addition, $G$ has an order-unit, then one can choose $R$
unital.\qed
\end{ET}

Direct limits of countable sequences of simplicial groups,
or, more generally, of arbitrary directed families of
simplicial groups, were characterized by
E.G. Effros, D.E. Handelman and C.-L. Shen as (countable)
dimension groups, see \cite{EHS80}. However, a very similar
result was proved four years earlier by P.A. Grillet
\cite{Gril76}, using a categorical result of R.T. Shannon
\cite{Shan74}. We refer to Section~\ref{S:DirDim} for
details. This characterization, together with Elliott's
Lemma, allows one to conclude that any countable dimension
group is isomorphic to $K_0$ of a locally matricial algebra
of countable dimension. That representation result was
extended by Handelman and the first author \cite{GoHa86} to
dimension groups of size $\aleph_1$.

\smallskip

The basic aim of this paper is to bring all these results
together. For instance, we prove, in
Theorem~\ref{T:ReprSemil}, the following analogue
of the Grillet and Effros-Handelman-Shen theorems:
\emph{every distributive semilattice is a direct limit of
finite Boolean semilattices}. This gives, in
Section~\ref{S:Berg}, a second proof of Bergman's Theorem.

A third proof of Bergman's Theorem, also given in
Section~\ref{S:Berg}, involves the relationship between
dimension groups and distributive semilattices. More
specifically, we prove in Theorem~\ref{T:LiftCtble} that
every countable distributive $0$-semilattice is isomorphic
to the maximal semilattice quotient of some countable
dimension group, and then we apply the Effros-Handelman-Shen
Theorem and Elliott's Lemma. The machinery that allows us
to conclude is, in fact, disseminated in the literature, and
it is recalled in Section~\ref{S:SCRR}.

A parallel to Bergman's Theorem, in which any distributive
algebraic lattice with countably many compact elements is
represented as the lattice of closed ideals of an
approximately finite-dimensional C*-algebra, is developed in
Section~\ref{S:Cstar}. As an application, we use this result
to provide a normal form for certain C*-algebras recently
classified by H. Lin in \cite{Lin97}.

Thus the present paper is, at the same time, a survey about
many intricately interwoven results in the theories of
dimension groups, semilattices, regular rings, C*-algebras,
and complemented modular lattices, which have been evolving
with various degrees of mutual independence for decades.

\section{Basic concepts}\label{S:BasConc}

We denote by $\omega$ the set of all natural numbers. A
natural number
$n$ is identified with the finite set
$\{0,1,\ldots,n-1\}$.

If $f\colon X\to Y$ is a map, $\ker(f)$, the \emph{kernel}
of $f$, denotes the equivalence relation associated with
$f$, that is,
 \[
 \ker(f)=\{\vv<u,v>\in X\times X\colon f(u)=f(v)\}.
 \]

\smallskip

We write \cm s additively, and we endow every \cm\ with its
\emph{algebraic preordering} $\leq$, defined by
 \[
 x\leq y\text{ if and only if there exists }z
 \text{ such that }x+z=y.
 \]
An \emph{ideal} (sometimes called an \emph{o-ideal}) of a
\cm\ $M$ is a nonempty subset $I$ of
$M$ such that for all $x$, $y\in M$, $x+y\in I$ if and only
if $x\in I$ and $y\in I$. (Note that this is a different
concept than the notion of `ideal' as used in semigroup theory.)
Write $\Id M$ for the set of ideals of
$M$, ordered by inclusion, and observe that $\Id M$ is a
complete lattice (with infima given by intersections).

The \emph{refinement property} is the semigroup-theoretical
axiom stating that for all positive integers $m$ and $n$,
all elements $a_i$ ($i<m$) and $b_j$ ($j<n$) of $M$ such
that $\sum_{i<m}a_i=\sum_{j<n}b_j$, there are elements
$c_{ij}$ ($i<m$, $j<n$) of $M$ such that
 \[
 a_i=\sum_{j<n}c_{ij}\quad\text{for all }i<m,
 \qquad\text{and}\qquad
 b_j=\sum_{i<m}c_{ij}\quad\text{for all }j<n.
 \]
A \emph{refinement monoid} (\emph{e.g.},
\cite{Dobb82}, \cite{Wehr92}) is a \cm\ which satisfies the
refinement property; equivalently, the condition above is
satisfied for $m=n=2$. It is to be noted that in
\cite{Dobb82}, every refinement monoid is, in addition,
required to satisfy the axiom $x+y=0\Rightarrow x=y=0$
(\emph{conicality}), while this is not the case for most other
authors (\emph{e.g.}, \cite{AGOMP}, \cite{Wehr92}).

\smallskip

A \emph{semilattice} is a commutative semigroup $S$ in which
every element $x$ is \emph{idempotent}, that is, $x+x=x$.
The algebraic preordering on $S$ is then an ordering, given
by $x\leq y$ if and only if $x+y=y$, hence all our semilattices
are join-semilattices.
We will usually denote by $\vee$, rather than $+$, the addition
of a semilattice. An \emph{ideal} (or \emph{order-ideal}) of $S$
is defined by the same axiom used to define an ideal of a
monoid. In order-theoretic terms, an ideal of $S$ is any nonempty
\emph{lower subset} $I$ (\emph{i.e.}, $(\forall x\in
S)(\forall y\in I) (x\le y \implies x\in I)$) which is
closed under $\vee$. A \emph{$0$-semilattice} is a semilattice
which is also a monoid, or, equivalently, a semilattice which
has a least element. Similarly, a \emph{$0$-lattice} is a
lattice with a least element.

An element $a$ of a lattice $L$ is \emph{compact} if, for
every subset $X$ of $L$ such that $\bigvee X$ exists, if
$a\leq\bigvee X$, then there exists a finite subset $Y$ of
$X$ such that $a\leq\bigvee Y$. Note that the set of compact
elements of $L$ forms a subsemilattice of $L$. A lattice~$L$
is \emph{algebraic} if $L$ is complete and every element of
$L$ is a supremum of compact elements.

If $S$ is a semilattice, denote by $\Id S$ the set of ideals
of $S$, ordered under inclusion. The \emph{canonical
embedding} from $S$ into $\Id S$ is defined by
 \[
 s\mapsto\dnw s=\{x\in S\colon x\leq s\}.
 \]
Observe that $\Id S$ is a lattice if and only if $S$ is
downward directed, and is a complete lattice if and only if
$S$ has a least element. In the latter case, $\Id S$ is an
algebraic lattice. Conversely, for every algebraic lattice
$L$, the set of all compact elements of $L$ is a
$0$-semilattice. The following classical result (\emph{cf.}
\cite[Theorem VIII.8]{Birk67}) expresses the
categorical equivalence between algebraic lattices and join
$0$-semilattices.

\begin{proposition}\label{AlgLattId}
Let $L$ be an algebraic
lattice, and let $S$ be the semilattice of all compact
elements of $L$. Then the correspondence
 \[
 x\mapsto\{s\in S\colon s\leq x\}
 \]
defines an isomorphism from $L$ onto $\Id S$.\qed
\end{proposition}

This can be extended without difficulty to define a
categorical equivalence between $0$-semilattices and
$\{\vee,0\}$-homomorphisms, and algebraic lattices with a
suitable notion of homomorphism.

A semilattice $S$ is \emph{distributive}
(see \cite[p.~117]{Grat71})
if for all $a$, $b_0$, $b_1$ in $S$
such that $a\leq b_0\vee b_1$, there are elements $a_0$ and
$a_1$ of $S$ such that $a=a_0\vee a_1$ and $a_i\leq b_i$ for
all $i<2$. This is equivalent to saying that $S$ is downward
directed and $\Id S$ is a distributive lattice,
\emph{cf.} \cite[Lemma 11.1(iii)]{Grat71} or
\cite[Lemma II.5.1]{Grat98}. Together with
Proposition~\ref{AlgLattId}, this shows that if $L$ is a
distributive algebraic lattice, then the semilattice of all
compact elements of $L$ is a distributive semilattice.

For every lattice $L$, we denote by $\Con L$ the lattice of
all congruences of $L$. It is a well known theorem of N. Funayama and
T. Nakayama (see \cite[Corollary 9.16]{Grat71} or
\cite[II.3]{Grat98}) that $\Con L$ is a distributive
algebraic lattice.

We denote by $\Conc L$ the
semilattice of
\emph{compact} congruences of $L$; by the previous
paragraph, $\Conc L$ is a distributive $0$-semilattice. The
elements of $\Conc L$ are exactly the finitely generated
congruences of $L$.

\smallskip

For every \poag\ $G$, we denote by $G^+$ the positive cone
of $G$, that is, the set of $x\in G$ such that $x\geq 0$. An
\emph{order-unit} of $G$ is any element
$u$ of $G^+$ such that for every $x\in G$, there exists a
positive integer $n$ such that $x\leq nu$. We put
$\NN=\ZZ^+\setminus\{0\}$.

Let $G$ and $H$ be \poag s. A \emph{positive homomorphism}
from $G$ to $H$ is a homomorphism of \poag s from $G$ to
$H$,  that is, a group homomorphism $f\colon G\to H$ such
that $f(G^+)\subseteq H^+$. We denote by $f^+$ the
restriction of $f$ from $G^+$ to $H^+$.

\smallskip

All the rings that we will consider are associative, but not
necessarily unital.

\section{Refinement monoids, dimension groups and\\
distributive semilattices}\label{S:RefDimDis}

Let $M$ be a \cm. There exists a least monoid congruence
$\asymp$ on $M$ such that $M/{\asymp}$ is a semilattice. It
is convenient to define $\asymp$ in terms of the preordering
$\propto$ defined by
 \[
 x\propto y\text{ if and only if }
 (\exists n\in\NN)(x\leq ny);
 \]
then, $x\asymp y$ if and only if $x\propto y$ and
$y\propto x$. The
\emph{maximal semilattice quotient} of $M$ is the natural
projection from $M$ to $M/{\asymp}$, often identified with
the semilattice $M/{\asymp}$ itself. We refer to \cite{ClPr}
for the details.

This defines a \emph{functor} from the category of \cm s,
with monoid homomorphisms, to the category of
$0$-semilattices, with $0$-semilattice homomorphisms. We
will denote this functor by $\nabla$. The proof of the
following lemma is straightforward.

\begin{lemma}\label{L:PiPresDirLim}
The functor $\nabla$ preserves direct limits.\qed
\end{lemma}

Now let us go to refinement monoids:

\begin{lemma}\label{L:RefRiesz}
Every refinement monoid $M$ satisfies the
\emph{Riesz decomposition property}, that is, for all
elements $a$, $b_0$ and $b_1$ of $M$ such that $a\leq
b_0+b_1$, there are $a_0\leq b_0$ and $a_1\leq b_1$ in $M$
such that $a=a_0+a_1$.\qed
\end{lemma}

For semilattices, it is well known (and also easy to verify
directly) that the converse of Lemma~\ref{L:RefRiesz} is
true:

\begin{lemma}\label{L:DistRef}
Let $S$ be a semilattice.
Then $S$ is distributive if and only if $S$ satisfies the
refinement property.\qed
\end{lemma}

We will be interested in the effect of $\nabla$ on refinement
monoids:

\begin{lemma}\label{L:PiRef}
Let $M$ be a \cm. If $M$ satisfies the Riesz decomposition
property, then $\nabla(M)$ is a distributive semilattice.
\end{lemma}

\begin{proof}
For every element $x$ of $M$, denote by $[x]$
the image of $x$ in $\nabla(M)$. Let $a$, $b_0$ and $b_1$ be
elements of $M$ such that $[a]\leq[b_0]\vee[b_1]$. By
definition, there exists $n\in\NN$ such that $a\leq
nb_0+nb_1$. Since $M$ satisfies the Riesz decomposition
property, there are
$a_0\leq nb_0$ and $a_1\leq nb_1$ such that $a=a_0+a_1$.
Therefore,
$[a]=[a_0]\vee[a_1]$, and $[a_i]\leq[b_i]$ for $i<2$.
\end{proof}

Say that a partially ordered set $\vv<P,\leq>$ satisfies the
\emph{interpolation property} if, for all $a_0$, $a_1$,
$b_0$ and $b_1$ in $P$ such that
$a_i\leq b_j$ for all $i$, $j<2$, there exists $x\in P$ such
that $a_i\leq x\leq b_j$ for all
$i$, $j<2$. An \emph{interpolation group} is a \poag\
satisfying the interpolation property.

\begin{lemma}
[{see \cite[Proposition 2.1]{Gpoag}}]
\label{L:IGChar}
Let $G$ be a \poag. Then $G$ is an
interpolation group if and only if its positive cone $G^+$
is a refinement monoid.\qed
\end{lemma}

Say that a \poag\ $G$ is \emph{directed} if it is directed as
a partially ordered set; equivalently, $G=G^++(-G^+)$. Say
that $G$ is \emph{unperforated} if for all $m\in\NN$ and all
$x\in G$, $mx\geq0$ implies that $x\geq0$. A
\emph{dimension group} is a directed, unperforated
interpolation group.

For example, define a \emph{dimension vector space} (over
$\QQ$) to be a directed interpolation group endowed with a
structure of vector space over the field $\QQ$ of the
rational numbers, for which multiplication by positive
rational scalars is order-preserving. Then it is obvious
that every dimension vector space is a dimension group
($mx\geq0$ implies $(1/m)mx\geq0$, that is,
$x\geq0$).

By Lemma~\ref{L:PiRef}, the maximal semilattice quotient of
the positive cone of a dimension group is a distributive
$0$-semilattice. The converse is an open problem (see
Problem~\ref{Pb:GrpSemil} in Section~\ref{S:Probs}).

In Sections~\ref{S:TempLex} and~\ref{S:Lifting}, we will
solve positively two particular cases of this problem: the
case where $S$ is a \emph{lattice}
(Theorem~\ref{T:LiftDLat}), and the case where $S$ is
\emph{countable} (Theorem~\ref{T:LiftCtble}). \smallskip

An \emph{ideal} of a \poag\ $G$ is a subgroup $I$ of $G$
which is both
\emph{directed} and
\emph{convex} with respect to the ordering on $G$, the
latter condition meaning that whenever
$x\leq y\leq z$ with
$x$, $z\in I$ and $y\in G$, then $y\in I$. We denote by $\Id G$
the set of ideals of $G$, ordered under inclusion; by
\cite[Corollary 1.10]{Gpoag}, $\Id G$ is a
complete lattice. Let $\Idc G$ denote the subsemilattice of
compact elements in $\Id G$. It is an easy exercise to see that
an ideal $I$ of $G$ lies in $\Idc G$ if and only if $I$ has an
order-unit (when $I$ is viewed as a \poag\ in its own right).

Similarly, for any \cm\ $M$ we write $\Idc M$ for the
semilattice of compact elements of $\Id M$, and we observe
that the members of $\Idc M$ are precisely those ideals of
$M$ which have order-units.

\begin{proposition}\label{IdStuff}
Let $M$ be a \cm\ and $G$ a \poag.

\begin{enumerate}
\item $\Id G\cong \Id G^+$.

\item $\Id M$ and $\Id G$ are algebraic lattices.

\item If $M$ satisfies the Riesz decomposition property,
then $\Id M$ is distributive and $\Idc M\cong \nabla(M)$.
Hence, $\Id M\cong \Id \nabla(M)$.

\item If $G$ is an interpolation group, then $\Id G$ is
distributive and $\Idc G\cong \nabla(G^+)$. Hence, $\Id
G\cong \Id \nabla(G^+)$.
\end{enumerate}
\end{proposition}

\begin{proof}
(i) Inverse isomorphisms are given as follows:
map each ideal $I$ of $G$ to $I\cap G^+$, and map each ideal
$J$ of $G^+$ to $J+(-J)$.

(ii) We already know that $\Id M$ and $\Id G$ are complete
lattices. Any ideal
$I$ of $M$ is the supremum of the \emph{principal} ideals
$\{y\in M\colon y\propto x\}$ for $x\in I$, and each of
these principal ideals is in $\Idc M$. This shows that $\Id
M$ is algebraic. One can argue similarly that $\Id G$ is
algebraic, or just apply part (i).

(iii) It follows directly from Riesz decomposition that for
any two ideals $I$ and $J$ of $M$, the sum $I+J$ is again an
ideal. Hence, finite suprema in $\Id M$ are given by sums. It
is clear that $(I+J)\cap K= (I\cap K)+(J\cap K)$ for all
$I$, $J$, $K\in\Id M$, and therefore $\Id M$ is distributive.

Observe that elements $x$, $y\in M$ satisfying $x\asymp y$
generate the same principal ideal of $M$. Hence, there is a
map $\theta\colon\nabla(M)\to \Idc M$ such that
$$\theta([x])= \{z\in M\colon z\propto x\}$$ for all $x\in
M$. Observe that $\theta([x])\subseteq \theta([y])$ if and
only if $x\propto y$, if and only if
$[x]\leq [y]$. Hence, $\theta$ is an order embedding. Any
ideal $I\in \Idc M$ has an order-unit, say $x$, and so $I=
\theta([x])$. Therefore $\theta$ is an order-isomorphism of
$\nabla(M)$ onto $\Idc M$, hence also a semilattice
isomorphism.

That $\Id M\cong \Id\nabla(M)$ now follows from
Proposition~\ref{AlgLattId}.

(iv) By Lemma~\ref{L:IGChar} and parts (i), (iii) above, we
have that $\Id G\cong \Id G^+$ is distributive
(\emph{cf.} \cite[Propositions 2.4, 2.5]{Gpoag})
and $\Idc G\cong \Idc G^+\cong \nabla(G^+)$. Now
$\Id G\cong\Id\nabla(G^+)$ by Proposition~\ref{AlgLattId}.
\end{proof}

\section{Direct limit representation of dimension groups;
Triangle Lemma}\label{S:DirDim}

Here we discuss the Effros-Handelman-Shen Theorem and
separate its proof into two parts: a ``Triangle Lemma''
concerning positive homomorphisms from simplicial groups to
dimension groups, and a ``direct limit representation
lemma'' which provides sufficient conditions for objects of
a quasivariety to be represented as direct limits of objects
from a given subclass. The latter lemma we prove in detail,
as it will yield our direct limit representation theorem for
distributive semilattices (Theorem~\ref{T:ReprSemil})
once we establish a suitable Triangle Lemma in that setting
(Corollary~\ref{C:TriLem}).

A \emph{simplicial group} is a \poag\ that is isomorphic to
some $\ZZ^n$, equipped with the direct product ordering, for
a nonnegative integer
$n$. Obviously, every simplicial group is a dimension group.
Conversely, it turns out that simplicial groups are
``building blocks" of dimension groups, \emph{via} direct limits.
The earliest result of this type is due to P.A. Grillet
\cite[Theorem 2.1]{Gril76}. Say that a \cm\ $S$ has the
\emph{strong Riesz interpolation property} (strong RIP) if for
every positive integer $n$ and for all elements $a$, $b$, $c$ and
$d$ of $S$, if $na+b=nc+d$, then there are elements $u$, $v$,
$w$ and $z$ of $S$ such that $a=u+v$, $b=nw+z$, $c=u+w$, and
$d=nv+z$.

\begin{theorem}[{see \cite[Theorem 2.1]{Gril76}}]
\label{T:Grillet}
Let $S$ be a \cm. Then the following are
equivalent:
\begin{enumerate}
\item $S$ is a direct limit of (finitely generated) free \cm
s.

\item $S$ is cancellative and $S$ satisfies the strong RIP.
\qed
\end{enumerate}
\end{theorem}

The passage from the strong RIP to the direct limit
representation is achieved by using a general categorical
result, due to R.T. Shannon \cite{Shan74}, which gives a
characterization of directed colimits of free objects in
algebraic categories.

\begin{remark}
Although the fact is absent from
\cite{Gril76}, it is not difficult, although not
trivial, to verify directly that a directed \poag\ $G$ is a
dimension group if and only if $G^+$ satisfies the strong RIP.
To establish the nontrivial implication, one starts by proving
directly that any dimension group $G$ satisfies
Proposition~3.23 of \cite{Gpoag}, that is, for
all $n\in\NN$ and all $a\in G$, the set of $x\in G$ such that
$a\leq nx$ is downward directed. This can be done by induction
on $n$; here is an outline of a proof.

Let $x_0$, $x_1\in G$ such that $a\leq nx_i$ for $i<2$, with
$n\geq 2$. For all $i$, $j<2$, we have
$na=a+(n-1)a\leq nx_i+(n-1)nx_j$ and so,
by $n$-unperforation, $a\le x_i+(n-1)x_j$. Apply interpolation
to the relations $a-x_i\leq (n-1)x_j$ to obtain $y\in G$ such
that $a-x_i\leq y\leq (n-1)x_j$ for all $i$, $j$. By the
induction hypothesis, there exists $z\in G$ such that
$y\leq(n-1)z$ and $z\leq x_j$ for
$j<2$. Since also $a-y\leq x_j$ for all $j$, another
interpolation yields $x\in G$ such that
$z\, ,a-y\leq x\leq x_0\, ,x_1$. Therefore
$a\leq x+y\leq x+(n-1)z\leq nx$, completing the induction step.

Now to prove the strong RIP, let
$n\in\NN$ and $a$, $b$, $c$ and $d$ in $G^+$ such that
$na+b=nc+d$. Put
$e=na-d=nc-b$, and note that $e\leq na,nc$. Thus there
exists $u$ such that
$e\leq nu$ and $u\leq a,c$. Put $v=a-u$, $w=c-u$, and
$z=d-nv$.

Therefore, this exercise is an easy proof that Grillet's
Theorem (Theorem~\ref{T:Grillet}) implies the later
Effros-Handelman-Shen Theorem (Theorem~\ref{T:EHS})
described in the next paragraph.

\end{remark}

The direct limit representation result for dimension groups
was proved by E.G. Effros, D.E. Handelman and C.-L.~Shen:

\begin{theorem}[{see \cite[Theorem 2.2]{EHS80}}]
\label{T:EHS}
A \poag\ is a direct limit of simplicial groups if and only
if it is a dimension group.\qed
\end{theorem}

A proof of this result is also presented in
\cite[Theorem 3.19]{Gpoag}. The hard core of the
proof consists in what we shall call the \emph{Triangle Lemma}:
For every simplicial group $S$, every dimension group
$G$ and every positive homomorphism $f\colon S\to G$, there
exist a simplicial group $T$ and positive homomorphisms
$\varphi\colon S\to T$ and $g\colon T\to G$ such that
$f=g\circ\varphi$ and
$\ker(f)=\ker(\varphi)$. Once this step is established, the
argument follows a general, categorical pattern. There are,
in fact, general categorical results which allow one to go
directly from the Triangle Lemma above to the direct limit
representation. For example, the main result of R.T. Shannon
\cite{Shan74} is quite short to state (modulo numerous
necessary definitions), but we did not find it convenient to
translate it, for example, to the language of \poag s for
the purpose of finding a shorter proof of the
Effros-Handelman-Shen Theorem. On the other hand, it seems
almost unavoidable that writing down the most general
categorical statement that leads from the Triangle Lemma to
the direct limit representation would involve a substantial
number of extremely unwieldy statements.

To solve this dilemma, we will put ourselves at a medium
level of generality, which will be sufficient to deal with
current first-order theories (such as \cm s, or
semilattices). While Shannon's result is stated in a
categorical context, we will choose a universal algebraic
context. This way, the reader can at least choose, according
to his affinities, between a categorical statement and a
universal algebraic statement.

We assume familiarity with only the very rudiments of
universal algebra, and we refer to \cite{Malc}
for the details. We will fix a language $\LL$ of
\emph{algebras}, that is, a first-order language with only
symbols of operations and constants (no relation symbols).
Say that a \emph{quasi-identity} is a first-order sentence
of the form
 \[
 (\forall\vec\vx)\bigl[
 \varphi(\vec\vx)\Rightarrow\psi(\vec\vx)\bigr],
 \]
where $\varphi$ is a finite (possibly empty)
conjunction of \emph{equations}
( = atomic formulas) and $\psi$ is an equation. A
\emph{quasivariety} (see \cite[Chapter V]{Malc})
is the class of models of a set of quasi-identities. It is
well-known that in any quasivariety $\VV$, there are
arbitrary colimits. In particular, for every set
$X$, there exists a free object of $\VV$ over $X$.

\begin{lemma}\label{L:TrtoLim}
Let $\VV$ be a quasivariety
of algebras of $\LL$, and let $M\in\VV$. Let $\FF$ be a
subclass of $\VV$ with the following properties:

\begin{enumerate}
\item For each $m\in M$, there exist $F\in\FF$ and a
homomorphism $f\colon F\to M$ such that $m\in f(F)$;

\item For each coproduct $F$ of finitely many elements
of $\FF$ and each homomorphism $f\colon F\to M$, there exist
$G\in\FF$ and homomorphisms $\varphi\colon F\to G$ and
$g\colon G\to M$ such that $f=g\circ\varphi$ and
$\ker(f)=\ker(\varphi)$.
\end{enumerate}

Then $M$ is a direct limit of objects from $\FF$.
\end{lemma}

\begin{proof}
We mimic the proof presented in
\cite[Theorem 3.19]{Gpoag}. Put $I=M\times\omega$. (This is just to
ensure that we base our indexing on an infinite set, to cover the
possibility that $M$ might be finite.) Put
$P=[I]^{<\omega}\setminus\{\varnothing\}$, the set of all
nonempty finite subsets of $I$, ordered under inclusion. We
construct inductively objects $F_p\in\FF$ and homomorphisms
$f_p\colon F_p\to M$ for $p\in P$, and transition
homomorphisms $f_{pq}\colon F_p\to F_q$ for $p\subset q$ in
$P$ (where $\subset$ denotes \emph{strict} inclusion). If
$p=\{\vv<a,n>\}$, where $\vv<a,n>\in M\times\omega$, choose, by
hypothesis (i), an
$F_p\in\FF$ and a homomorphism $f_p\colon F_p\to M$ such
that $a\in f_p(F_p)$.

Now the induction step. Suppose that $p\in P$ has at least
two elements, and suppose that we have constructed objects
$F_q\in\FF$ for $q\subset p$ in $P$, homomorphisms
$f_q\colon F_q\to M$ for $q\subset p$ in $P$, and
$f_{qr}\colon F_q\to F_r$ for $q\subset r\subset p$ in $P$,
satisfying the following conditions:

\begin{enumerate}
\item If $p_0\subset p_1\subset p_2\subset p$ in $P$, then
$f_{p_0p_2}=f_{p_1p_2}\circ f_{p_0p_1}$.

\item If $p_0\subset p_1\subset p$ in $P$, then
$f_{p_0}=f_{p_1}\circ f_{p_0p_1}$.

\item If $p_0\subset p_1\subset p$ in $P$, then
$\ker(f_{p_0})=\ker(f_{p_0p_1})$.
\end{enumerate}

Put $F^p=\coprod_{q\subset p}F_q$, where $\amalg$ denotes
the coproduct in
$\VV$. For all $q\subset p$ in $P$, denote by $e_{qp}$ the
canonical homomorphism from $F_q$ to $F^p$. By the universal
property of the coproduct, there exists a unique
homomorphism $f^p\colon F^p\to M$ such that $f^p\circ
e_{qp}=f_q$ for all $q\subset p$ in $P$. By assumption,
there exist an object $F_p\in\FF$ and homomorphisms
$\varphi_p\colon F^p\to F_p$ and $f_p\colon F_p\to M$ such
that $f_p\circ\varphi_p=f^p$ and $\ker(f^p)=\ker(\varphi_p)$.
For all $q\subset p$ in $P$, define
$f_{qp}=\varphi_p\circ e_{qp}$. The construction may be
described by the commutative diagram below:

\begin{picture}(100,70)(-50,-15)

\put(10,5){\vector(2,1){70}}
\put(90,5){\vector(0,1){35}}
\put(10,0){\vector(1,0){70}}
\put(100,0){\vector(1,0){70}}
\put(170,5){\vector(-2,1){70}}

\put(92,22){\makebox(0,0)[l]{$f^p$}}
\put(40,-5){\makebox(0,0)[t]{$e_{qp}$}}
\put(40,35){\makebox(0,0)[t]{$f_q$}}
\put(130,-5){\makebox(0,0)[t]{$\varphi_p$}}
\put(135,35){\makebox(0,0)[t]{$f_p$}}

\put(0,0){\makebox(0,0){$F_q$}}
\put(90,0){\makebox(0,0){$F^p$}}
\put(90,45){\makebox(0,0){$M$}}
\put(180,0){\makebox(0,0){$F_p$}}

\end{picture}

We verify points (i) to (iii) listed above for the larger set
of all $q\in P$ such that $q\subseteq p$.

(i) It suffices to verify that, for
$p_0\subset p_1\subset p$, we have
$f_{p_0p}=f_{p_1p}\circ f_{p_0p_1}$, that is,
$\varphi_p\circ e_{p_0p}=
\varphi_p\circ e_{p_1p}\circ f_{p_0p_1}$. Since
$\ker(\varphi_p)=\ker(f^p)$, it suffices to prove that
$f^p\circ e_{p_0p}=f^p\circ e_{p_1p}\circ f_{p_0p_1}$, that
is, $f_{p_0}=f_{p_1}\circ f_{p_0p_1}$, which is indeed the
case by the induction hypothesis (ii).

(ii) It suffices to verify that, for $q\subset p$ in $P$, we
have $f_q=f_p\circ f_{qp}$. This is a direct calculation:
 \[
 f_p\circ f_{qp}=f_p\circ\varphi_p\circ e_{qp}=
 f^p\circ e_{qp}=f_q.
 \]

(iii) It suffices to verify that, for $q\subset p$ in $P$, we
have $\ker(f_q)=\ker(f_{qp})$. Let $x$, $y\in F_q$. Then
$f_{qp}(x)=f_{qp}(y)$ if and only if
$\varphi_p\circ e_{qp}(x)=\varphi_p\circ e_{qp}(y)$. Since
$\ker(\varphi_p)=\ker(f^p)$, this is equivalent to $f^p\circ
e_{qp}(x)=f^p\circ e_{qp}(y)$, that is, $f_q(x)=f_q(y)$.

Therefore, we have constructed a direct system
 \[
 \mathcal{S}=\vv<\vv<F_p,f_{pq}>\colon
 p\subset q\text{ in }P>,
 \]
and homomorphisms $f_p\colon F_p\to M$ such that
$f_p=f_q\circ f_{pq}$ and $\ker(f_p)=\ker(f_{pq})$ for all
$p\subset q$ in $P$. Further, for each $a\in M$ we have,
for $p=\{\vv<a,0>\}$, that $p\in P$ and $a\in f_p(F_p)$.
Now if $S$, together with
limiting maps $\eta_p\colon F_p\to S$, is the direct limit of the system
$\mathcal{S}$ in $\VV$, there exists a unique homomorphism $f\colon S\to
M$ such that $f\circ \eta_p= f_p$ for all $p\in P$, and $f$ is
surjective. To see that $f$ is injective, let $\vv<s_0,s_1>\in \ker(f)$.
Then there exist $p\in P$ and $x_0$, $x_1\in F_p$ such that
$\eta_p(x_i)=s_i$ for $i<2$, and $\vv<x_0,x_1>\in \ker(f_p)$. Since $I$
is infinite, there exists $q\in P$ such that
$p\subset q$, and $\vv<x_0,x_1>\in \ker(f_{pq})$ by construction, whence
$s_0= (\eta_q\circ f_{pq})(x_0)= (\eta_q\circ f_{pq})(x_1)= s_1$.
Therefore $f$ is an isomorphism.
\end{proof}

\subsection*{Example}
In the language consisting of a binary
operation symbol $+$ and a constant symbol $0$, one can
consider the quasivariety of \emph{\cm s}. Finitely
generated free \cm s are exactly the positive cones of
simplicial groups. The Triangle Lemma in this context is a
reformulation of the corresponding Triangle Lemma for
\poag s (Shen's condition), see for example
\cite[Proposition 3.16]{Gpoag}. It is to be
noted that Lemma~\ref{L:TrtoLim} cannot be directly applied to
\poag s, because of the binary
\emph{relation} symbol $\leq$. However, this is easily
finessed here by considering the positive cones instead of
the full ordered groups.

We will see another application of Lemma~\ref{L:TrtoLim} in
Section~\ref{S:RepSemil}, in the case of the variety of
\emph{semilattices}.

\section{Temperate powers of $\QQ$}\label{S:TempLex}

The purpose of this section is to demonstrate that
Problem~\ref{Pb:GrpSemil} (see Section~\ref{S:Probs}) has a
positive solution for distributive \emph{$0$-lattices}.
Moreover, the construction developed here will allow us, in the
following section, to demonstrate a positive solution to
Problem~\ref{Pb:GrpSemil} for countable distributive
semilattices.

Throughout this section, we shall fix a set $X$ and a
sublattice $D$ of the powerset lattice $\mathcal{P}(X)$,
such that $\varnothing\in D$. Define
$\BB(D)$ to be the generalized Boolean subalgebra of
$\mathcal{P}(X)$ generated by $D$. Equivalently, the
elements of $\BB(D)$ are finite unions of the form
 \[
 \bigcup_{i<2n}\bigl(a_{2i}\setminus a_{2i+1}\bigr),
 \]
where $\vv<a_0,a_1,\ldots,a_{2n}>$ is a finite decreasing
sequence of elements of $D$ (see \cite[II.4]{Grat98}).

Further, let $\EE{D}$ be the set of all functions $f\colon
X\to\QQ$ with finite range such that $f$ is
\emph{measurable} with respect to the generalized Boolean
algebra $\BB(D)$, that is, $f^{-1}\{r\}$ belongs to $\BB(D)$
for every nonzero $r\in\QQ$.

\begin{lemma}\label{L:EDjVectSp}
The set $\EE{D}$ is a subalgebra of the
$\QQ$-algebra $\QQ^X$. Furthermore, for all $f$,
$g\in\EE{D}$, the map $(f:g)$, defined componentwise by
 \begin{equation*}
 (f:g)(x)=
   \begin{cases}
     f(x)/g(x)&\text{if }g(x)\ne 0,\\
     0&\text{if }g(x)=0
   \end{cases}
 \end{equation*}
belongs to $\EE{D}$.
\end{lemma}

\begin{proof}
It is obvious that $\EE{D}$ is closed under
multiplication by rational scalars.

For all $f$, $g\in\EE{D}$, both $f$ and
$g$ have finite range, thus so does $f+g$. Furthermore, for
all $r\in\QQ$, we have
 \[
 (f+g)^{-1}\{r\}=
 \bigcup\bigl\{ f^{-1}\{s\}\cap g^{-1}\{t\}\colon
 s\in\rng(f),\ t\in\rng(g)\text{ and }s+t=r\bigr\},
 \]
where $\rng(f)$ denotes the range of $f$, and thus
$(f+g)^{-1}\{r\}\in\BB(D)$. Hence,
$f+g\in\EE{D}$. Similarly, the product $fg$ and the element
$(f:g)$ belong to $\EE{D}$.
\end{proof}

For every element $f$ of $\EE{D}$, note that the
\emph{support} of $f$,
 \[
 \supp(f)=\{x\in X\colon f(x)\ne 0\},
 \]
belongs to $\BB(D)$.

\begin{notation}
Let $\EEp{D}$ be the set of all functions
$f\in\EE{D}$ such that $f(x)\geq0$ for all $x\in X$, and
$\supp(f)\in D$.
\end{notation}

\begin{proposition}\label{P:EDjpovs}
$\EEp{D}$ is the positive cone of a structure of dimension
vector space on $\EE{D}$.
\end{proposition}

\begin{proof}
It is easy to verify that $\EEp{D}$ is the
positive cone of a structure of partially ordered vector
space on $\EE{D}$: one has to verify that
$\EEp{D}$ is an additive submonoid of $\EEp{D}$, closed
under multiplication by positive rational numbers, and that
$\EEp{D}\cap(-\EEp{D})=\{0\}$; this is straightforward.

Denote by $\leq$ the pointwise ordering of $\EE{D}$, and by
$\leq^+$ the ordering of $\EE{D}$ with positive cone
$\EEp{D}$.

Every element $f\in\EE{D}$ is majorized (for $\leq$) by some
$n\cdot\chi_A$, where $n\in\NN$ and $A\in D$ (here $\chi_Y$
denotes the characteristic function of a subset $Y$ of $X$).
Therefore, the support of $(n+1)\cdot\chi_A-f$ is equal to $A$,
so that $f\leq^+(n+1)\cdot\chi_A$. Hence, the partial ordering
$\leq^+$ is directed.

It remains to verify interpolation. It is convenient to use
Lemma~\ref{L:IGChar}, that is, to verify that $\EEp{D}$
satisfies the refinement property.

Thus let $f_0$, $f_1$, $g_0$ and $g_1$ be elements of
$\EEp{D}$ such that
$f_0+f_1=g_0+g_1$. Put $h=f_0+f_1$. For all $i$, $j<2$, put
(with the notation of Lemma~\ref{L:EDjVectSp})
 \[
 h_{ij}=(f_ig_j:h).
 \]

By Lemma~\ref{L:EDjVectSp}, $h_{ij}$ belongs to $\EE{D}$. It
is obvious that
$0\leq h_{ij}$. To prove that $0\leq^+h_{ij}$, note that
$\supp(h_{ij})=\supp(f_i)\cap\supp(g_j)$. But $D$ is closed
under finite intersection, whence $\supp(h_{ij})\in D$. Thus
$h_{ij}\in\EEp{D}$. Finally, it is obvious that
$f_i=h_{i0}+h_{i1}$ and
$g_j=h_{0j}+h_{1j}$ for $i$, $j<2$.
\end{proof}

In the sequel, we shall identify in notation $\EE{D}$ with
the dimension vector space
$\vv<\EE{D},+,0,\EEp{D}>$, and we will call it the
\emph{temperate power} of
$\QQ$ by $D$.

\begin{lemma}\label{L:PropEDj}
Let $f$ and $g$ be two
elements of $\EEp{D}$. Then the following are equivalent:
\begin{itemize}
\item[\rm(i)] There exists $n\in\NN$ such that $f\leq^+ng$.

\item[\rm(ii)] There exists $n\in\NN$ such that $f\leq ng$.

\item[\rm(iii)] $\supp(f)\subseteq\supp(g)$.

\end{itemize}
\end{lemma}

In particular, for $f$, $g\in\EEp{D}$, there is no ambiguity
on the notation
$f\propto g$, whether $\leq$ or $\leq^+$ is used to order
the vector space
$\EE{D}$.

\begin{proof}
(i)$\Rightarrow$(ii) and (ii)$\Rightarrow$(iii) are trivial.

Assume (iii). Since $(f:g)$ has finite range, it is
majorized by some positive integer $n$. Let $x\in X$; we
prove that $f(x)\leq ng(x)$. This is trivial when
$f(x)=0$. If $x\in\supp(f)$, that is, $f(x)>0$, then, by
assumption, $g(x)>0$, thus $(f:g)(x)=f(x)/g(x)$; but
$n\geq(f:g)(x)$, so that $f(x)\leq ng(x)$. It follows easily
that the support of $(n+1)g-f$ is equal to the support of
$g$; whence $f\leq^+(n+1)g$.
\end{proof}

By putting together Proposition~\ref{P:EDjpovs} and
Lemma~\ref{L:PropEDj}, one obtains the following result:

\begin{theorem}\label{T:LiftDLat}
For every distributive $0$-lattice $D$, there exists a
dimension vector space $E$ such that $\nabla(E^+)$ is
isomorphic to $D$ as a semilattice.
\end{theorem}

\begin{proof}
By Stone's Theorem (see \cite[Corollary II.1.21]{Grat98}),
there exists a set $X$ such that $D$ embeds into
$\mathcal{P}(X)$. Since $D$ has a zero, the embedding can be
arranged in such a way that its range includes
$\varnothing$. Thus, without loss of generality, we may assume
that $D$ is a sublattice of
$\mathcal{P}(X)$ containing $\varnothing$. Put $E=\EE{D}$.
By Proposition~\ref{P:EDjpovs}, $E$ is a dimension vector
space. By Lemma~\ref{L:PropEDj}, the maximal semilattice
quotient of $E^+$ is isomorphic to $D$ (\emph{via} the
support map).
\end{proof}

\section{Lifting countable distributive semilattices to
dimension groups}\label{S:Lifting}

In this section, we shall see how an easy application of the
results of Section~\ref{S:TempLex} yields a solution of
Problem~\ref{Pb:GrpSemil} in the case of \emph{countable}
semilattices.

For every partially ordered set $P$, denote by $\HH(P)$ the
distributive lattice of all lower subsets of $P$ (that is,
the subsets $X$ of $P$ such that if $p\leq x$ and $x\in X$,
then $p\in X$). Put
$\GG{P}=\EE{\HH(P)}$, and $\GGp{P}=\EEp{\HH(P)}$.
We will call
$\GG{P}$ (with positive cone $\GGp{P}$) the
\emph{temperate power} of $\QQ$ by $P$.

In case $P$ is finite, one can give a direct description of
the dimension vector space $\GG{P}$, since the generalized
Boolean algebra $\BB(\HH(P))$ just equals
$\mathcal{P}(P)$ in this case. The underlying space of
$\GG{P}$ is $\QQ^P$, and
$\GGp{P}$ consists of those functions $u\colon P\to \QQ^+$
whose support belongs to
$\HH(P)$.

By Lemma~\ref{L:PropEDj}, one can define an
\emph{isomorphism} $\iota_P\colon\nabla(\GGp{P})\to\HH(P)$,
by the formula
 \[
 \iota_P([x])=\supp(x),\qquad\text{for all }x\in\GGp{P}.
 \]

\begin{lemma}\label{L:LiftMaps}
Let $P$ and $Q$ be two finite partially ordered sets, and let
$f\colon\HH(P)\to\HH(Q)$ be a $0$-semilattice homomorphism.
Then there exists a positive homomorphism
$g\colon\GG{P}\to\GG{Q}$ such that
$\iota_Q\circ\nabla(g^+)=f\circ\iota_P$.
\end{lemma}

The last condition of the statement above means that the
following diagram commutes:
 \[
 \begin{CD}
 \nabla(\GGp{P})@>\iota_P>>\HH(P)\\ @V\nabla(g^+)VV @VVfV\\
 \nabla(\GGp{Q})@>\iota_Q>>\HH(Q)
 \end{CD}
 \]

\begin{proof}
Denote by $\vv<\dot p\colon p\in P>$ the
canonical basis of $\QQ^P$, where $\dot p=\chi_{\{p\}}$, and by
$\vv<\dot q\colon q\in Q>$ the canonical basis of
$\QQ^Q$. Let
$g$ be the unique linear map from $\QQ^P$ to $\QQ^Q$ defined
by the formula
 \[
 g(\dot p)=\sum_{q\in f(\dnw p)}\dot q
 \]
(recall the notation $\dnw p=\{x\in P\colon x\leq p\}$).
Let $x\in\GGp{P}$, written as $x=\sum_{p\in P}x_p\dot p$,
where all $x_p$ are elements of $\QQ^+$. Then we have
 \[
 g(x)=\sum_{p\in P}
 \left(x_p\sum_{q\in f(\dnw p)}\dot q\right)=
 \sum_{q\in Q}y_q\dot q,
 \]
where we put $D_q=\{p\in P\colon q\in f(\dnw p)\}$ and
$y_q=\sum_{p\in D_q}x_p$ for all $q\in Q$. It is obvious
that all $y_q$ belong to $\QQ^+$.

Put $U=\supp(x)$; by assumption, $U$ belongs to $\HH(P)$. If
$q\in f(U)$, then, since $U=\bigcup_{p\in U}\dnw p$ and
since $f$ is a $0$-semilattice homomorphism, there exists
$p\in U$ such that $q\in f(\dnw p)$, that is, $p\in D_q$.
Since $p\in U$, we have $x_p>0$, whence $y_q>0$. Conversely,
if $q\notin f(U)$, then, for all $p\in D_q$, we have
$p\notin U$ and thus $x_p=0$; hence, $y_q=0$. This shows
that $\supp(g(x)) =f(U)\in \HH(Q)$.

It follows that $g$ is a positive homomorphism, and that,
for all $x\in\GGp{P}$, we have
 \[
 \supp(g(x))=f(\supp(x)).
 \]
Hence, $g$ satisfies the required condition.
\end{proof}

By using Pudl\'ak's Lemma (see the Introduction), we can now
conclude:

\begin{theorem}\label{T:LiftCtble}
Every countable distributive $0$-semilattice $S$ is
isomorphic to the maximal semilattice quotient of the
positive cone of some countable dimension vector space $E$.
If, in addition, $S$ is bounded, then $E$ has an order-unit.
\end{theorem}

\begin{proof}
Let $S$ be a countable distributive
$0$-semilattice. By Pudl\'ak's result, one can write $S$ as
a countable, increasing union
$S=\bigcup_{n\in\omega}S_n$, where all the $S_n$ are finite
distributive subsemilattices of $S$, containing $0$. Then
each $S_n$ is a distributive lattice. Denote by $P_n$ the
set of all (nonzero) join-irreducible elements of $S_n$,
ordered by the restriction of the ordering of $S_n$, and by
$\tau_n$ the natural isomorphism from $\HH(P_n)$ onto $S_n$.
Put $f_n=\tau_{n+1}^{-1}|_{S_n}\circ\tau_n$. By
Lemma~\ref{L:LiftMaps}, there exists a positive homomorphism
$g_n\colon\GG{P_n}\to\GG{P_{n+1}}$ such that
$\iota_{P_{n+1}}\circ\nabla(g_n^+)=f_n\circ\iota_{P_n}$. The
information can be partly visualized in the following
commutative diagram:
 \[
 \begin{CD}
 \nabla(\GGp{P_n})@>\iota_{P_n}>>\HH(P_n)@>\tau_n>> S_n\\
 @V\nabla(g_n^+)VV @VVf_nV @VV\subseteq V\\
 \nabla(\GGp{P_{n+1}})@>\iota_{P_{n+1}}>>\HH(P_{n+1})
 @>\tau_{n+1}>>S_{n+1}\\
 \end{CD}
 \]
Consider the direct system $\mathcal{S}$ of partially
ordered $\QQ$-vector spaces whose objects are the
$\GG{P_n}$, for $n\in\omega$, and whose morphisms are the
maps $g_{n-1}\circ\cdots\circ g_m$, for $m\leq n$. By
Lemma~\ref{L:PiPresDirLim}, if $E$ denotes the direct limit
of $\mathcal{S}$, then $\nabla(E^+)$ is isomorphic to the
direct limit of the $S_n$ with the inclusion maps, that is,
to $S$. Since all the $\GG{P_n}$ are dimension vector spaces
(by Proposition~\ref{P:EDjpovs}), so is
$E$. It is clear that $E$ is countable.

Finally, suppose that $S$ is bounded. Then $\nabla(E^+)$ has
a largest element, call it $1$. Let $u\in E^+$ be an element
whose $\asymp$-class is $1$. Hence, all elements $x\in E^+$
satisfy $x\propto u$, and so $u$ is an order-unit of the
monoid $E^+$. Since $E$ is directed, $u$ must also be an
order-unit for $E$.
\end{proof}

\section{Boolean direct limit representation of distributive
semilattices}\label{S:RepSemil}

The Triangle Lemma for distributive semilattices can be
proved in a very similar fashion as the corresponding result
for dimension groups (\emph{i.e.}, \cite[Proposition 3.16]{Gpoag}).
However, we present here a different proof, that shows at the
same time a stronger property of distributive semilattices
(Proposition~\ref{P:FinInj}). Furthermore, this proof is
specific to semilattices, \emph{e.g.}, the analogue of
Proposition~\ref{P:FinInj} for dimension groups and
simplicial groups does \emph{not} hold.

\begin{lemma}\label{L:ab+x}
Let $S$ be a distributive semilattice. Let
$n\in\omega$ and let $a$, $b$, $c_i$ ($i<n$) be elements of
$S$ such that $a\leq b\vee c_i$ for all $i<n$. Then there
exists $x\in S$ such that $a\leq b\vee x$ and $x\leq c_i$
for all $i<n$.
\end{lemma}

\begin{proof}
It suffices to prove the lemma for $n=2$.
Since $S$ is distributive, there are $b_i\leq b$ and
$d_i\leq c_i$ such that $a=b_i\vee d_i$ (for all $i<2$).
Since $d_1\leq a\leq b_0\vee d_0$, there are, further,
$b_2\leq b_0$ and $x\leq d_0$ such that $d_1=b_2\vee x$.
Therefore $x\leq d_i\leq c_i$ for all
$i<2$, and
 \begin{equation*}
 a=b_1\vee b_2\vee x\leq b\vee x.\tag*{\qed}
 \end{equation*}
\renewcommand{\qed}{}
\end{proof}

\begin{lemma}\label{L:EqSysDis}
Let $S$ be a distributive semilattice. Let $m$,
$n\in\omega$, and let $a_i$, $b_i$ ($i<m$) and $c_j$ ($j<n$)
be elements of $S$ such that $a_i\leq b_i\vee c_j$ for all
$i<m$ and $j<n$. Then there exists $x\in S$ such that
 \[
 (\forall i<m)(a_i\leq b_i\vee x)\quad\text{and}\quad
 (\forall j<n)(x\leq c_j).
 \]
\end{lemma}

\begin{proof}
This is an immediate consequence of \cite[Lemma 1.5]{ShWe94}.
However, we present here a self-contained proof.

By Lemma~\ref{L:ab+x}, for all
$i<m$, there exists $x_i\in S$ such that $a_i\leq b_i\vee
x_i$ and
$x_i\leq c_j$ for all $j<n$. Then
$x=\bigvee_{i<m}x_i$ satisfies the required conditions.
\end{proof}

\begin{proposition}[Finite injectivity for distributive
semilattices]
\label{P:FinInj}
Let $S$ be a distributive semilattice, and
let $A$ be a subsemilattice of a finite semilattice $B$.
Then every semilattice homomorphism from $A$ to $S$ extends
to a semilattice homomorphism from $B$ to $S$.
\end{proposition}

\begin{proof}
Let $f$ be a homomorphism from $A$ to $S$.

We consider first the case where there exists $b\in
B\setminus A$ such that $B$ is generated by $A\cup\{b\}$.
Therefore,
 \begin{equation}\label{Eq:BAb}
 B=A\cup\{b\}\cup\{x\vee b\colon x\in A\}.
 \end{equation}
Let $\{\vv<x_i,y_i>\colon
i<m\}$ list all the pairs $\vv<x,y>$ of elements of $A$ such
that $x\leq y\vee b$, and let $\{z_j\colon j<n\}$ list all
elements $z$ of $A$ such that $b\leq z$. For all $i<m$ and
all $j<n$, we have $x_i\leq y_i\vee z_j$, and thus
$f(x_i)\leq f(y_i)\vee f(z_j)$. Therefore, by
Lemma~\ref{L:EqSysDis}, there exists $\alpha\in S$ such that
 \begin{equation*}
  \begin{aligned}
  f(x_i)&\leq f(y_i)\vee\alpha&&\qquad\quad
  \text{ for all }i<m,\\
  \alpha&\leq f(z_j)&&\qquad\quad\text{ for all }j<n.
  \end{aligned}
 \end{equation*}
Let $\{w_k\colon k<r\}$ list all elements $w$ of $A$ such
that $w\le b$. Then
$f(w_k)\le f(z_j)$ for all $k<r$ and $j<n$. Set
$\beta=\alpha\vee\bigvee_{k<r}f(w_k)$
($\beta$ is defined as being equal to $\alpha$ if $r=0$),
and observe that
 \begin{equation}\label{Eq:finj2}
  \begin{aligned}
   f(x_i)&\leq f(y_i)\vee\beta&&\qquad\quad
   \text{ for all }i<m,\\
   f(w_k)&\leq\beta&&\qquad\quad\text{ for all }k<r,\\
   \beta&\leq f(z_j)&&\qquad\quad\text{ for all }j<n.
  \end{aligned}
 \end{equation}
It follows from (\ref{Eq:finj2}) that
 \begin{equation}\label{Eq:finj3}
   \begin{aligned}
   f(x)\vee\beta &=f(y)\vee\beta &&\qquad\quad
   \text{for all } x,y\in A
   \text{ such that } x\vee b =y\vee b\\
   f(x)\vee\beta &=f(z)&&\qquad\quad\text{for all }
   x,z\in A \text{ such that }x\vee b =z\\
   f(x)\vee\beta&=\beta &&\qquad
   \quad\text{for all } x\in A
   \text{ such that\ } x\vee b =b.
   \end{aligned}
 \end{equation}
By (\ref{Eq:BAb}) and (\ref{Eq:finj3}),
$f$ extends to a well-defined map $g\colon B\to S$ such that
$g(b)=\beta$ and $g(x\vee b)= f(x)\vee\beta$ for all $x\in A$.
Since $f$ is a homomorphism, it follows easily that $g$ is a
homomorphism.

In the general case, there exists a finite chain of
subsemilattices
 \[
 A=B_0\subset B_1\subset\cdots\subset B_k=B
 \]
such that each $B_i$ $(i>0)$ is generated by
$B_{i-1}\cup\{b_i\}$ for some $b_i\in B_i\setminus B_{i-1}$.
Thus, we conclude by an easy induction argument.
\end{proof}

It is to be noted that Proposition~\ref{P:FinInj} is also an
immediate consequence of Pudl\'ak's Lemma (see the
Introduction) and the \emph{injectivity} of every finite
distributive semilattice in the class of semilattices. The
latter result follows immediately from
\cite[Theorem 3.11]{Wehr92}, but it can also be proved
directly. Moreover, our proof here is self-contained.

A finite semilattice is \emph{Boolean} if it is isomorphic
to $\tw^n$ for some $n\in\omega$, where $\tw$ is the two
element semilattice.

\begin{lemma}[folklore]
\label{L:EmbFin}
Every finite semilattice ($0$-semilattice) has a
(zero-pre\-serv\-ing) embedding into a finite Boolean lattice.
\end{lemma}

\begin{proof}
If $S$ is a finite semilattice, let
$B=\mathcal{P}(S)$ be the powerset semilattice of $S$, and
embed $S$ into $B$ \emph{via} the map $j\colon S\to B$
defined by the rule
 \begin{equation*}
 j(s)=\{x\in S\colon s\nleq x\}.\tag*{\qed}
 \end{equation*}
\renewcommand{\qed}{}
\end{proof}

A better embedding (from the computational viewpoint) can be
obtained by replacing $\mathcal{P}(S)$ by
$\mathcal{P}(P)$, where $P$ denotes the set of meet-irreducible
elements of $S$ (here meet-irreducibility means with respect to
whatever meets might exist); the map
$j$ is defined similarly. We can now prove the Triangle Lemma for
distributive semilattices:

\begin{corollary}\label{C:TriLem}
Let $S$ be a distributive semilattice. Let $A$ be a finite
semilattice, and let $f$ be a homomorphism from $A$
to $S$. Then there exist a finite Boolean semilattice $B$ and
homomorphisms $\varphi\colon A\to B$ and $g\colon B\to S$
such that $f=g\circ\varphi$ and $\ker(f)=\ker(\varphi)$.
\end{corollary}

\begin{proof}
Put $A'=A/\ker(f)$, and denote by
$\pi\colon A\twoheadrightarrow A'$ the quotient map. There
exists a unique homomorphism $f'\colon A'\to S$ such that
$f=f'\circ\pi$. By Lemma~\ref{L:EmbFin}, there exists an
embedding $j$ from $A'$ into some finite Boolean semilattice
$B$. By Proposition~\ref{P:FinInj}, there exists a
homomorphism $g\colon B\to S$ such that $f'=g\circ j$. Put
$\varphi=j\circ\pi$. The situation can be described by the
following commutative diagram:

\begin{picture}(100,70)(-50,-15)

\put(10,5){\vector(2,1){70}}
\put(90,5){\vector(0,1){35}}
\put(10,0){\vector(1,0){70}}
\put(100,0){\vector(1,0){70}}
\put(170,5){\vector(-2,1){70}}

\put(92,22){\makebox(0,0)[l]{$f'$}}
\put(40,-5){\makebox(0,0)[t]{$\pi$}}
\put(40,35){\makebox(0,0)[t]{$f$}}
\put(130,-5){\makebox(0,0)[t]{$j$}}
\put(135,35){\makebox(0,0)[t]{$g$}}

\put(0,0){\makebox(0,0){$A$}}
\put(90,0){\makebox(0,0){$A'$}}
\put(90,45){\makebox(0,0){$S$}}
\put(180,0){\makebox(0,0){$B$}}

\end{picture}

We obtain the following:
 \[
 f=f'\circ\pi=g\circ\varphi.
 \]
Furthermore, $j$ is one-to-one, and thus
$\ker(\varphi)=\ker(\pi)=\ker(f)$.
\end{proof}

We can now deduce a general representation result for
distributive semilattices:

\begin{theorem}\label{T:ReprSemil}
Every distributive semilattice is a direct limit of
finite Boolean semilattices and semilattice homomorphisms.
\end{theorem}

\begin{proof}
We consider the first-order language consisting of one binary
operation symbol $\vee$, the variety $\VV$ of semilattices, and
the subclass $\FF$ of finite Boolean semilattices. Since
the class of finite (not necessarily Boolean) semilattices is
closed under finite coproducts (because every finitely generated
semilattice is finite), the assumption (ii) of
Lemma~\ref{L:TrtoLim} is, by Corollary~\ref{C:TriLem},
satisfied. Since the assumption (i) of Lemma~\ref{L:TrtoLim}
is trivially satisfied, the theorem follows.
\end{proof}

Say that a partially ordered set is \emph{bounded} if it has
a least and a greatest element, which we denote by $0$ and
$1$.

\begin{corollary}\label{C:ReprSemil01}\hfill
\begin{enumerate}
\item Every distributive $0$-semilattice is a direct limit of
finite Boolean semilattices and $0$-preserving
semilattice homomorphisms.

\item Every bounded distributive semilattice is a direct limit
of finite Boolean semilattices and $0,1$-preserving
semilattice homomorphisms.

\end{enumerate}
\end{corollary}

\begin{proof}
We prove, for example, (i). The proof for (ii) is similar.
Let $S$ be a distributive $0$-semilattice.
By Theorem~\ref{T:ReprSemil}, $S$ is a direct limit of a
direct system
 \[
 \mathcal{S}=\vv<\vv<S_i,f_{ij}>\colon i\leq j\text{ in }I>
 \]
where $I$ is a directed set, the $S_i$ are finite Boolean
semilattices and the $f_{ij}$ are semilattice
homomorphisms, with respect to limiting homomorphisms
$f_i\colon S_i\to S$.

Without loss of generality, $I$ has a least element, denoted
by $0$, and $0_S=f_0(0_{S_0})$. For all $i\in I$, put
$0_i=f_{0i}(0_{S_0})$ and
$T_i=\{x\in S_i\colon 0_i\leq x\}$. Then $T_i$ is a finite
Boolean semilattice, and $f_{ij}$ maps $T_i$ to $T_j$ for
$i\leq j$. Furthermore, the least element of $T_i$ is
$0_i$, and $f_{ij}(0_i)=0_j$ for $i\leq j$. Thus each $f_{ij}$
restricts to a
$0$-preserving semilattice homomorphism
$g_{ij}\colon T_i\to T_j$. Finally,
$S$ is the direct limit of the system
 \begin{equation*}
 \mathcal{T}=\vv<\vv<T_i,g_{ij}>\colon i\leq j\text{ in }I>.
 \tag*{\qed}
 \end{equation*}
\renewcommand{\qed}{}
\end{proof}

\begin{example}
Consider the three element chain $S=\{0,1,2\}$, viewed as a
bounded join-semilattice. Although
Corollary~\ref{C:ReprSemil01} allows us to express $S$ as a
direct limit of finite Boolean semilattices, the result is
puzzling, because $S$ itself, although finite, is not Boolean.

Here is an explicit description of $S$ as a direct limit of
finite Boolean semilattices. Consider the $0$-semilattice
homomorphism $r\colon\tw^2\to\tw^2$ defined by $r(a)=a$ and
$r(b)=a\vee b$, where $a$ and $b$ are the two atoms of $\tw^2$.
It is not difficult to verify that $S$ is the direct limit of
the sequence
 \begin{equation*}
 \begin{CD}
 \tw^2 @>r>> \tw^2 @>r>> \tw^2 @>r>> \cdots,
 \end{CD}
 \end{equation*}
with the limiting $0$-semilattice homomorphism
$f\colon\tw^2\to S$ defined by $f(a)=1$ and $f(b)=2$.
\end{example}

\section{Regular rings and the functors $V$,
$\Vtil$}\label{S:SCRR}

We recall the definition and some basic facts about regular
rings, their idempotents, and their ideal lattices.

For every ring $R$, denote by $\LL(R)$ the semilattice of
all finitely generated right ideals of $R$, ordered by
inclusion. A ring $R$ is (von Neumann)
\emph{regular} if for all $x\in R$, there exists $y\in R$
such that $xyx=x$.

A $0$-lattice $L$ is \emph{sectionally complemented} if for
all elements $a\leq b$ of $L$, there exists a sectional
complement of $a$ in $b$, that is, an element $x$ of $L$ such
that $a\wedge x=0$ and $a\vee x=b$.

\begin{proposition}\label{P:RegSC}
If $R$ is a regular ring, then $\LL(R)$ is a sectionally
complemented modular lattice.
\end{proposition}

\begin{proof}
This was first proved by von~Neumann in the unital case
\cite[Theorem 2]{vonN36}; his argument easily
extends to the non-unital situation, as noted in
\cite[3.2]{FrHa56}.
\end{proof}

Let $R$ be a ring. For all $n\in\NN$, embed the ring
$M_n(R)$ of all $n\times n$ square matrices over $R$ into
$M_{n+1}(R)$, \emph{via} the map
 \[
 x\mapsto
 \begin{pmatrix}
   x&0\\
   0&0
 \end{pmatrix} .
 \]
Furthermore, denote by $M_\infty(R)$ the direct limit of
the system
 \[
 R\to M_2(R)\to M_3(R)\to M_4(R)\to\cdots.
 \]
Define an
equivalence relation
$\sim$ on the set of all idempotent elements of
$M_\infty(R)$ by
 \[
 e\sim f\Longleftrightarrow (\exists x,y\in M_\infty(R))
 (xy=e\text{ and }yx=f).
 \]
Equivalently, $e\sim f$ if and only if $e\cdot
M_\infty(R)\cong f\cdot M_\infty(R)$ as right
$M_\infty(R)$-modules. For every idempotent $e$ of
$M_\infty(R)$, denote by $[e]$ the $\sim$-equivalence class
of $a$, and put
 \[
 V(R)=\{\,[e]\colon e\in R,\ e^2=e\,\}.
 \]
There is a well-defined addition on $V(R)$ given by
 \begin{align*}
 [e]+[f]=[e\oplus f]\\
 \intertext{where}
   e\oplus f=
   \begin{pmatrix}
     e&0\\
     0&f
   \end{pmatrix},
 \end{align*}
and $\vv<V(R),+,[0]>$ is a \cm. Now, $V$
extends to a \emph{functor} from the category of rings to
the category of \cm s. It is well known (and also easy to
see) that this functor preserves finite direct products and
direct limits.

Since we shall often work with the maximal semilattice
quotients of the monoids
$V(R)$, let us introduce the notation $\Vtil$ for the
composition of $V$ with the functor $\nabla$ (see
Section~\ref{S:RefDimDis}). Thus $\Vtil$ is a functor from
the category of rings to the category of semilattices, and it
preserves direct limits and finite direct products. Given a
ring~$R$, write $|e|$ for the
$\asymp$-class of $[e]\in V(R)$, where $e$ is any idempotent
in $M_\infty(R)$.

\begin{proposition}\label{P:RegRef}
Let $R$ be a regular ring. Then $V(R)$ is a
\crm, and $\Vtil(R)$ is a distributive semilattice.
\end{proposition}

\begin{proof}
That $V(R)$ is a \crm\ follows from Theorem~2.8 in
\cite{Gvnrr} (the fact that $R$ does not
necessarily have a unit does not affect the proof). By
Lemma~\ref{L:PiRef},
$\Vtil(R)$ is a distributive semilattice.
\end{proof}

It is well known that for any ring $R$, the lattice $\Id R$
of (two-sided) ideals of $R$ is algebraic. The semilattice
$\Idc R$ of compact elements of $\Id R$ consists of the
finitely generated ideals of
$R$, that is, all two-sided ideals of $R$ of the form
$\sum_{i<n}Rx_iR$, where
$n\in\NN$ and all $x_i$ belong to $R$. Note that by
Proposition~\ref{AlgLattId}, $\Id R\cong \Id(\Idc R)$.

\begin{proposition}\label{P:ManyIso}
Let $R$ be a regular
ring. Then all three semilattices $\Conc\LL(R)$, $\Idc R$
and $\Vtil(R)$ are pairwise isomorphic. Furthermore, they
are distributive $0$-semi\-lat\-tices. Moreover, the
lattices $\Id R$, $\Id V(R)$, and $\Id\Vtil(R)$ are pairwise
isomorphic, and these are distributive algebraic lattices.
\end{proposition}

\begin{proof}
The semilattice isomorphisms follow from
\cite[Corollary 4.4 and Proposition 4.6]{Wehr2}. Again, the
fact that $R$ does not necessarily have a unit does not
affect the proofs. Then by Proposition~\ref{P:RegRef}, these
semilattices are distributive. (Recall that more generally, the
semilattice of all compact congruences of any lattice is
distributive, see Section~\ref{S:BasConc}).

Now $\Id R\cong \Id(\Idc R)\cong \Id\Vtil(R)\cong \Id V(R)$.
By Proposition~\ref{IdStuff}, $\Id\Vtil(R)$ is distributive,
and the proposition is proved. (One can also prove directly
that $\Id R$ is distributive; this is well known and easy.)
\end{proof}

As a byproduct of the proof of Proposition~\ref{P:ManyIso},
we have the following:

\begin{proposition}\label{P:IdVRIR}
Let $R$ be a regular ring.
\begin{enumerate}
\item Let $J$ be an ideal of $V(R)$. Then there exists a
two-sided ideal $I$ of $R$ such that $J\cong V(I)$. Namely,
$I$ is the ideal of $R$ generated by all idempotents $e\in R$
for which $[e]\in J$.

\item Let $J$ be an ideal of $\Idc R$. Then there exists a
two-sided ideal $I$ of $R$ such that $J\cong\Idc I$. Namely,
$I$ is the sum of all those ideals of
$R$ which are members of $J$.\qed
\end{enumerate}
\end{proposition}

Note that every two-sided ideal of a regular ring is itself
regular (see \cite[Lemma 1.3]{Gvnrr}).

Another consequence of Proposition~\ref{P:ManyIso} is the
following:

\begin{corollary}\label{C:RegComp}
Let $D$ be a distributive algebraic lattice. If there exists
a regular ring $R$ such that $D\cong\Id R$, then there exists
a sectionally complemented modular lattice $L$ such that
$D\cong\Con L$.
\end{corollary}

Translated into the language of semilattices, this gives the
following: \emph{Let $S$ be a distributive $0$-semilattice.
If there exists a regular ring $R$ such that $S\cong\Idc R$,
then there exists a sectionally complemented modular lattice
$L$ such that $S\cong\Conc L$.}

\begin{proof}
Propositions~\ref{AlgLattId}, \ref{P:RegSC},
and \ref{P:ManyIso}.
\end{proof}

\section{Bergman's Theorem}\label{S:Berg}

We are now ready to develop our two new proofs of Bergman's
Theorem. Let us first recall some basic definitions. Let $K$
be a field. A
\emph{matricial algebra} over $K$ is a finite direct product
of the form
 \[
 \prod_{i<k}M_{n_i}(K),
 \]
where $k$ and the $n_i$ are natural numbers. A
\emph{locally matricial} algebra over $K$ is a direct limit
of matricial algebras over $K$ and
$K$-algebra homomorphisms. Note that we do not require the
ring homomorphisms to preserve the ring units. Note also
that locally matricial algebras are very special cases of
regular rings. Countable dimensional locally matricial
algebras are sometimes called \emph{ultramatricial}, see
\cite{Gvnrr}.

Observe that if $R$ is a matricial algebra, then
$V(R)\cong (\ZZ^+)^n$ for some positive integer $n$
(see \cite[Lemma 15.22]{Gvnrr} for an analogous
result with the same proof). In particular, $V(R)$ is then
cancellative ($x+a=x+b$ implies $a=b$). Since the functor $V$
preserves direct limits, $V(R)$ is also cancellative for any
locally matricial algebra $R$. Thus $V(R)\cong K_0(R)^+$ for any
such $R$, since $K_0(R)$ is constructed as the universal
enveloping group of $V(R)$. We shall use this observation to
translate results from the literature, stated in the
language of $K_0$, into $V(R)$-form.

Elliott's Lemma (see the Introduction) together with the
countable case of the Effros-Handelman-Shen Theorem
(Theorem~\ref{T:EHS}), which implies that every countable
dimension group is the direct limit of a countable sequence
of simplicial groups, yields the following result:

\begin{theorem}[{see \cite[2nd.~Ed., p.~376]{Gvnrr}}]
\label{C:ElEHS}
Let $G$ be a countable dimension group, and let $K$ be a
field. Then there exists a locally matricial $K$-algebra $R$
of countable dimension such that $V(R)\cong G^+$. If, in
addition, $G$ has an order-unit, then one can choose $R$
unital.\qed
\end{theorem}

In this section, we will illustrate the interdependency of
various parts of this paper, by giving two proofs of
Bergman's Theorem (stated in the Introduction).

\begin{proof}[First Proof of Bergman's Theorem]
By Proposition~\ref{AlgLattId}, it suffices to solve the
following problem. We fix a countable distributive
$0$-semilattice $S$ and a field
$K$; we must find a locally matricial $K$-algebra $R$ of
countable dimension such that $\Idc R\cong S$. In view of
Proposition~\ref{P:ManyIso}, this is the same as to arrange
for $\Vtil(R)\cong S$. Further, if $S$ is bounded, we must
find a unital such $R$.

By Theorem~\ref{T:LiftCtble}, there exists a countable
dimension vector space
$E$ such that $\nabla(E^+)\cong S$. By
Theorem~\ref{C:ElEHS}, there exists a locally matricial
$K$-algebra $R$ of countable dimension such that $V(R)\cong
E^+$; therefore $\Vtil(R)\cong\nabla(E^+) \cong S$.

In addition, if $S$ is bounded, then $E$ has an order-unit,
and thus, by Theorem~\ref{C:ElEHS}, one can choose $R$
unital.
\end{proof}

\begin{proof}[Second Proof of Bergman's Theorem]
This proof does not use the results of Elliott, or Grillet,
Effros, Handelman and Shen. In fact, it uses nothing more
than the countable case of Corollary~\ref{C:ReprSemil01}.

As in the first proof, we fix a countable distributive
$0$-semilattice $S$ and a field $K$, and we find a locally
matricial $K$-algebra $R$ of countable dimension such that
$\Vtil(R)\cong S$.

According to Corollary~\ref{C:ReprSemil01}, we may assume that
$S$ is the direct limit of a sequence
\begin{equation}\label{CD:Seq}
\begin{CD}
\tw^{n_1} @>f_1>> \tw^{n_2} @>f_2>> \tw^{n_3} @>f_3>> \cdots
\end{CD}
\end{equation}
in the category of $0$-semilattices.

Set $R_1= K^{n_1}$ (the direct product of $n_1$ copies of
$K$), and observe that
$\vtil(R_1)\cong\tw^{n_1}$. More precisely, if $p_1,\dots,
p_{n_1}$ are the primitive central idempotents in $R_1$
(that is, the atoms of the finite Boolean algebra of central
idempotents of $R_1$), then $|p_1|,\dots,|p_{n_1}|$ are
distinct atoms which generate $\vtil(R_1)$. Hence, if
$a_1,\dots, a_{n_1}$ are the distinct atoms in $\tw^{n_1}$,
there exists an isomorphism $g_1\colon\vtil(R_1) \rightarrow
\tw^{n_1}$ such that
$g_1(|p_i|)=a_i$ for $i=1, \dots, n_1$.

Let $b_1,\dots, b_{n_2}$ be the distinct atoms in
$\tw^{n_2}$. There are integers $s_{ij}\in \{0,1\}$ such
that $f_1(a_i)=s_{i1}b_1 +\dots+ s_{i,n_2}b_{n_2}$ for all
$i$. Choose a positive integer $t(j)\ge s_{1j} +\dots+
s_{n_1,j}$ for each $j$, and set
 \[
 R_2= M_{t(1)}(K)\times\cdots\times M_{t(n_2)}(K).
 \]
Let $\phi_1\colon R_1\rightarrow R_2$ be the block diagonal
$K$-algebra homomorphism with multiplicities $s_{ij}$, that
is, each component map
$R_1\rightarrow M_{t(j)}(K)$ is given by
 \[
 (\alpha_1, \dots, \alpha_{n_1})\mapsto
 \diag(\overbrace{\alpha_1}^{s_{1j}},
 \overbrace{\alpha_2}^{s_{2j}}, \dots,
 \overbrace{\alpha_n}^{s_{n_1,j}}, 0,\dots,0),
 \]
where the notation
$\overbrace{\alpha_i}^{s_{ij}}$ means that $\alpha_i$
appears if $s_{ij}=1$ but not if $s_{ij}=0$. Let
 \[
 q_1=(I_{t(1)},0,\dots,0),\quad
 q_2=(0,I_{t(2)},0,\dots,0),\quad
 \dots, \quad q_{n_2}=(0,\dots,0,I_{t(n_2)})
 \]
be the primitive central idempotents in $R_2$. Then there
exists an isomorphism
$g_2\colon\vtil(R_2) \rightarrow \tw^{n_2}$ such that
$g_2(|q_j|)= b_j$ for all
$j$, and we observe that $g_2\vtil(\phi_1)= f_1g_1$.

Continuing in the same manner, we obtain a sequence
 \begin{equation}\label{CD:SeqR}
   \begin{CD}
   R_1 @>\phi_1>> R_2 @>\phi_2>> R_3 @>\phi_3>>\cdots
   \end{CD}
 \end{equation}
of matricial $K$-algebras and $K$-algebra
homomorphisms together with $0$-semilattice isomorphisms
$g_i\colon\vtil(R_i)\rightarrow\tw^{n_i}$ such that the
following diagram commutes:
 \[
  \begin{CD}
  \vtil(R_1) @>\vtil(\phi_1)>> \vtil(R_2) @>\vtil(\phi_2)>>
  \vtil(R_3) @>\vtil(\phi_3)>> \cdots\\
  @V g_1 VV @V g_2 VV @V g_3 VV \\
  \tw^{n_1} @>f_1>> \tw^{n_2} @>f_2>> \tw^{n_3} @>f_3>> \cdots
  \end{CD}
 \]
Therefore, if $R$ is the direct limit of the sequence
(\ref{CD:SeqR}), then we have $\vtil(R)\cong S$ as desired.

It remains to modify the proof for the case that $S$ has a
greatest element, say~$1$. As before, we express $S$ as the direct
limit of the sequence (\ref{CD:Seq}); in view of
Corollary~\ref{C:ReprSemil01}, we may now assume that the
maps $f_i$ preserve greatest elements. Thus
$f_i(1_i)=1_{i+1}$ for all $i$, where
$1_i$ denotes the greatest element of $\tw^{n_i}$ (the sum
of all the atoms).

Define $R_1$ as before, and note that $g_1$ maps $|1_{R_1}|$
to $1_1$.

Let $b_1,\dots, b_{n_2}$ and the $s_{ij}$ be as before.
Since
 \[
 1_2= f_1(1_1)=
 \sum_{i=1}^{n_1} f_1(a_i)=
 \sum_{i=1}^{n_1} \sum_{j=1}^{n_2}s_{ij}b_j,
 \]
we must have
$\sum_{i=1}^{n_1} s_{ij}\ne 0$ for all $j$, and so we can
choose $t(j)= s_{1j} +\dots+ s_{n_1,j}$. Now if $R_2$ and
$\phi_1$ are defined as before, $\phi_1$ is a unital
homomorphism.

Continuing as before, we can obtain a sequence
(\ref{CD:SeqR}) in which all the homomorphisms $\phi_i$ are
unital, and therefore $R$ is a unital algebra.
\end{proof}

\section{Ideal lattices in C*-algebras}\label{S:Cstar}

In this section, we use the methods of the previous section
to derive an analogue
of Bergman's Theorem for C*-algebras.
This result, in turn, has an interesting application to a
class of C*-algebras $A$ which have been classified by H. Lin
\cite{Lin97} in terms of the invariants
$V(A)$, which for this particular class are actually
distributive semilattices.

Throughout, we deal only with \emph{complex} C*-algebras.
Recall that the natural morphisms in the category of
C*-algebras are \emph{*-homomorphisms} ($\CC$-algebra
homomorphisms which preserve the involution *), since such
maps are automatically contractions with respect to
C*-algebra norms (see, \emph{e.g.},
\cite[Theorem 2.1.7]{Mur90}). Every finite-dimensional
C*-algebra has the form
 \[
 \prod_{i<k}M_{n_i}(\CC),
 \]
where the matrix algebras $M_{n_i}(\CC)$ are equipped
with the conjugate transpose involution and the operator
norm (see, \emph{e.g.}, \cite[Theorem III.1.1]{Dav96} or
\cite[Theorem 6.3.8]{Mur90}). A C*-algebra is said to be
\emph{AF} (for ``approximately finite-dimensional'') if it is
isomorphic to a direct limit (in the category of C*-algebras)
of a countable sequence of finite-dimensional C*-algebras and
*-homomorphisms.

We shall need the fact that the functor $V$ commutes with
C*-algebra direct limits (\emph{i.e.}, norm-completions of
*-algebra direct limits), see
\cite[5.2.4]{Blac86}. Since $V$ of any finite-dimensional
C*-algebra is obviously cancellative, it follows that $V(A)$ is
also cancellative for all AF C*-algebras $A$. Hence,
$V(A)\cong K_0(A)^+$ when $A$ is AF. It also follows that
$K_0$ of any AF C*-algebra is a countable dimension group,
see \cite[Theorem IV.3.3]{Dav96}.

In the category of C*-algebras, kernels correspond to
\emph{closed} ideals (ideals closed in the norm topology).
Thus, the natural ideal lattice to study is the lattice
$\Idbar A$ of closed ideals of a C*-algebra $A$. Such a
lattice is algebraic: infima are given by intersections,
suprema are given by closures of sums, and the compact
elements are the finitely generated closed ideals. (For an
ideal to be finitely generated in the context of closed
ideals means that it is the closure of some ideal which is
finitely generated in the usual sense.) It is known that the
lattice of closed ideals of an AF C*-algebra is distributive.
The C*-algebra analogue
of Bergman's Theorem can be stated as
follows:

\begin{theorem}\label{T:CstarBergman}
Let $L$ be a distributive algebraic lattice with only countably
many compact elements. Then $L$ is isomorphic to the lattice of
closed ideals in some AF C*-algebra $A$. If, in addition,
the greatest element of $L$ is compact, then $A$ can be
chosen to be unital.
\end{theorem}

\begin{proof}[First proof of Theorem~\ref{T:CstarBergman}]
By Proposition~\ref{AlgLattId}, there is a countable
distributive $0$-semilattice
$S$ such that $\Id S\cong L$. If $A$ is an AF C*-algebra,
let $\Idcbar A$ denote the semilattice of finitely generated
closed ideals of $A$. This semilattice consists precisely of
the compact elements of $\Idbar A$, and hence $\Id(\Idcbar
A)\cong \Idbar A$. Thus, it suffices to find an AF C*-algebra
$A$ such that $\Idcbar A\cong S$.

For any AF C*-algebra $A$, the lattice
$\Idbar A$ is isomorphic to $\Id K_0(A)$,
see \cite[Proposition IV.5.1]{Dav96}, and,
consequently, $\Idcbar A$ is isomorphic to $\Idc K_0(A)$. By
Proposition~\ref{IdStuff},
$\Idc K_0(A)\cong \nabla(K_0(A)^+)\cong \Vtil(A)$, and hence
$\Idcbar A\cong
\Vtil(A)$. Thus, to find an AF C*-algebra $A$ with $\Idbar
A\cong L$ is the same as to find an $A$ with $\Vtil(A)\cong
S$.

By Theorem~\ref{T:LiftCtble}, there exists a countable
dimension vector space $E$ such that
$\nabla(E^+)\cong S$. By the Effros-Handelman-Shen Theorem
and the C*-algebra analogue
of Elliott's Lemma, see
\cite[Theorem IV.7.3]{Dav96}, there exists an AF C*-algebra $A$
such that $V(A)\cong E^+$. Therefore $\Vtil(A)\cong S$,
as desired. In addition, if $S$ is bounded, then $E$ has
an order-unit, and then $A$ can be chosen to be unital.
\end{proof}

\begin{proof}[Second proof of Theorem~\ref{T:CstarBergman}]
As above, we just need to find an AF C*-algebra $A$ such
that $\Vtil(A)$ is isomorphic to a given countable
distributive $0$-semilattice $S$.

The construction in our second proof of Bergman's Theorem
yields a sequence (\ref{CD:SeqR}) of matricial
$\CC$-algebras and $\CC$-algebra homomorphisms such that the
direct limit of $\Vtil$ of (\ref{CD:SeqR}) is isomorphic to
$S$. Each $R_i$ can be viewed as a finite-dimensional
C*-algebra. Observe that the block diagonal maps $\phi_i$
are *-homomorphisms. Hence, the C*-algebra direct limit of
the sequence (\ref{CD:SeqR}) is an AF C*-algebra, say $A$.
Since the functor $V$ commutes with C*-algebra direct
limits, we therefore have $\Vtil(A)\cong S$, as desired.
\end{proof}

The result of Theorem~\ref{T:CstarBergman} can be extended
to other classes of C*-algebras by a simple tensor product
argument. For the basic theory of C*-tensor products and the
fundamental concept of \emph{nuclearity}, we refer the
reader to \cite[Chapter 6]{Mur90}. We shall need
the fact that all AF C*-algebras are nuclear, see
\cite[Theorem 6.3.11]{Mur90}. Since all the
C*-tensor products we consider will have at least one nuclear
factor, the C*-tensor products will be unique, and we will just
denote them by
$\otimes$.

The following lemma is well known among the cognoscenti, but
we have been unable to locate a reference in the literature,
and so we outline a proof here. We thank Bruce Blackadar for
this argument.

A C*-algebra $B$ is said to be \emph{simple} provided $B$ is
nonzero and the only closed ideals of $B$ are $0$ and $B$.

\begin{lemma}\label{L:TensorId}
Let $A$ and $B$ be C*-algebras, at least one of which is
nuclear. If $B$ is simple and unital, then
$\Idbar A\cong\Idbar(A\otimes B)$, \emph{via} the map
$I\mapsto I\otimes B$.
\end{lemma}

\begin{proof}
The rule $I\mapsto I\otimes B$ defines an order-preserving map
$\theta$ from $\Idbar A$ to $\Idbar(A\otimes B)$. Since $B$
is unital, there is a *-homomorphism $\phi\colon A\to
A\otimes B$ given by the rule $\phi(a)= a\otimes1$, and the
set map $\phi^{-1}$ induces an order-preserving map
$\varphi$ from $\Idbar(A\otimes B)$ to $\Idbar A$. Clearly
$\varphi\theta$ is the identity on $\Idbar A$. Thus, to prove
that $\theta$ is a lattice isomorphism, it suffices to show
that $\theta$ is surjective.

Let $J\in\Idbar(A\otimes B)$, set $I=\phi^{-1}(J)$, and
consider the algebraic (\emph{i.e.}, uncompleted) tensor
products $R= A\otimes_{\text{alg}} B$ and $S=
(A/I)\otimes_{\text{alg}} B$. Note that $K= J\cap R$ is an
ideal of $R$ such that $K\cap (A\otimes 1)= I\otimes 1$.
Since $B$ is simple and unital, its center is a field as well
as a C*-algebra, so the center of $B$ is $\CC\cdot1$.
Consequently, $K=I\otimes_{\text{alg}} B$ (see, \emph{e.g.},
\cite[Theorem V.6.1]{Jac64}).

Thus the composition of the inclusion map $R\to A\otimes B$
with the quotient map $A\otimes B\to (A\otimes B)/J$ induces
a *-algebra \emph{embedding}
$\psi\colon S\to (A\otimes B)/J$. The composition of $\psi$
with the quotient norm on $(A\otimes B)/J$ then defines a
C*-norm, call it $\Vert\cdot\Vert_\psi$, on $S$. (It is a
norm, rather than just a seminorm, because $\psi$ is
injective.) By, \emph{e.g.}, \cite[Theorem T.6.21]{Wegg93},
$\Vert\cdot\Vert_\psi$ is a cross norm on $S$. Because of
our nuclearity assumption, $\Vert\cdot\Vert_\psi$ is the
unique C*-cross norm on $S$, and so the completion of $S$
with respect to $\Vert\cdot\Vert_\psi$ yields the C*-tensor
product $(A/I)\otimes B$. On the other hand, $\psi$ is an
isometry and the image of $\psi$ is dense in $(A\otimes
B)/J$. Hence, $\psi$ induces a *-isomorphism of
$(A/I)\otimes B$ onto $(A\otimes B)/J$. It follows that the
kernel of the induced map $A\otimes B\to (A/I)\otimes B$ is
precisely $J$, and therefore $J= I\otimes B$, as desired.
\end{proof}

\begin{corollary}\label{C:TensorLatt}
Let $B$ be a simple, unital C*-algebra, and let $L$ be a
distributive algebraic lattice with only countably many compact
elements. Then there exists an AF C*-algebra $A$ such that
$L\cong\Idbar(A\otimes B)$. If, in addition, the greatest
element of $L$ is compact, then $A$ can be chosen to be unital.
\end{corollary}

\begin{proof}
Theorem~\ref{T:CstarBergman} and Lemma~\ref{L:TensorId}.
\end{proof}

We will apply the above corollary with a special choice of
$B$ which will ensure that $V(A\otimes B)$ is a distributive
semilattice. This is the \emph{Cuntz algebra} $\OO$, defined
as the unital C*-algebra generated by elements $s_1$ and
$s_2$ satisfying the relations
 \[
 s^*_1s_1= s^*_2s_2= s_1s^*_1+s_2s^*_2= 1.
 \]
It is known
that $\OO$ is simple (see, \emph{e.g.},
\cite[Corollary V.4.7]{Dav96}), that all nonzero projections in
$\OO$ are equivalent, see \cite[Corollary 3.12]{Cntz81}, and
that $M_n(\OO)\cong\OO$ for all $n\in\omega$, see
\cite{PaSa79}. In particular, it follows that $V(\OO)\cong\tw$.

In \cite{Lin97}, Lin classified a class $\AAA$ of
C*-algebras which have trivial K-theory, that is, the groups
$K_0$ and $K_1$ are both trivial for the algebras in $\AAA$.
We shall not give the precise definition of $\AAA$ here, but
just recall that $\AAA$ contains $\OO$ and is closed under
the following operations: hereditary C*-subalgebras,
quotients, tensor products with AF C*-algebras, countable
direct limits, finite tensor products, and extensions,
see \cite[Theorem 3.14]{Lin97}. For any $A\in\AAA$, the
monoid $V(A)$ is a countable distributive semilattice
(\emph{cf.} \cite[Proposition 3.4]{Lin97} and
\cite[Theorem 1.1]{Zha90}, or
\cite[Theorem 7.2 and Corollary 1.3]{AGOMP}). Further,
$\Idbar A\cong \Id V(A)$ for any $A\in\AAA$ (use
\cite[Proposition 3.4 and Corollary 3.11]{Lin97} to see that
$A$ has real rank zero, then use the argument of
\cite[Theorem 2.3]{Zha90}). Lin proved that the algebras in
$\AAA$ are classified up to isomorphism by the semilattices
$V(\cdot)$ together with elements corresponding to
approximate identities, see \cite[Theorem 3.13]{Lin97}.
In particular, unital algebras $A$, $B\in\AAA$ are isomorphic
if and only if $V(A)\cong V(B)$ [ibid].

Taking $B=\OO$ in Corollary~\ref{C:TensorLatt}, we see that
any distributive algebraic lattice with only countably many
compact elements can be represented as the lattice of closed
ideals of a C*-algebra $A\otimes\OO$ where $A$ is AF. Such tensor
products are in Lin's class $\AAA$, and we can use the above
information to see that all the unital algebras in $\AAA$
must have this form.

\begin{theorem}\label{T:LinUnital}
Each unital C*-algebra in Lin's class $\AAA$ has the form
$A\otimes\OO$ for some unital AF C*-algebra $A$.
\end{theorem}

\begin{proof}
Let $C\in\AAA$ be unital; then $V(C)$ is a bounded, countable,
distributive $0$-semilattice. The lattice $L=\Id V(C)$ is a
distributive algebraic lattice whose semilattice of compact
elements is isomorphic to $V(C)$ and thus is countable. Further,
the greatest element of $L$ is compact. By
Corollary~\ref{C:TensorLatt}, there exists a unital AF
C*-algebra $A$ such that
$L\cong\Idbar(A\otimes \OO)$. Hence,
 \[
 V(C)\cong \Idcbar(A\otimes\OO)\cong
 \Idc V(A\otimes\OO)\cong V(A\otimes\OO)/{\asymp}=
 V(A\otimes\OO).
 \]
Therefore we conclude from Lin's classification theorem,
see \cite[Theorem 3.13]{Lin97} that $C\cong A\otimes\OO$.
\end{proof}

\section{Open problems}\label{S:Probs}

The first circulated versions of the present paper 
generated some amount of work, which led to solutions to
most of the original open problems. The first one of these
open problems was the following.

\begin{problem}[Lifting distributive semilattices to
dimension groups]\label{Pb:GrpSemil}
Let $S$ be a distributive $0$-semilattice. Does there exist a
dimension group $G$ such that the maximal semilattice quotient
of $G^+$ (that is, $\nabla(G^+)$) is isomorphic to $S$?
\end{problem}

By Theorems~\ref{T:LiftDLat} and~\ref{T:LiftCtble},
Problem~\ref{Pb:GrpSemil} has positive solutions in the
\emph{lattice} case and in the \emph{countable} case. By
results of the first author and D.E. Handelman, see
\cite[Proposition 1.3]{Gartnoeth} and
\cite[second Corollary]{Hndlmn}, Problem~\ref{Pb:GrpSemil}
has a positive solution in case $S$ is totally ordered,
or---more generally---if every element of $S$ is a (finite)
join of join-irreducible elements of $S$.

On the other hand, P. R\r u\v zi\v cka
solved Problem~\ref{Pb:GrpSemil} negatively in \cite{Ruzi2},
for a semilattice $S$ of size $\aleph_2$. The $\aleph_1$ case
is still open:

\begin{all}{Problem~\ref{Pb:GrpSemil}$'$}
Let $S$ be a distributive $0$-semilattice of size $\aleph_1$.
Does there exist a dimension group $G$ such that
$\nabla(G^+)$ is isomorphic to $S$?
\end{all}

The regular ring version of Problem~\ref{Pb:GrpSemil}$'$ was
the following:

\begin{problem}[Lifting distributive semilattices to
regular rings]\label{Pb:aleph1Ring}
Let $S$ be a dis\-trib\-u\-tive $0$-semilattice. If
$|S|\leq\aleph_1$, does there exists a regular ring $R$ such
that $\Idc R$ is isomorphic to $S$?
\end{problem}

It is to be noted that the size $\aleph_1$ in the statement
of Problem \ref{Pb:aleph1Ring} is \emph{optimal}:
in~\cite{Wehr2}, the second author proved that there exists a
distributive $0$-semilattice $S$ of size $\aleph_2$ that
cannot be isomorphic to $\Idc R$ for any regular ring $R$. One
positive case of Problem~\ref{Pb:aleph1Ring} is that in which
$S$ is bounded and every element of $S$ is a finite join of
join-irreducible elements (no cardinality restriction on $S$
is needed). This follows from work of G.M. Bergman
\cite[\S\S2--4]{Berg86} extending the result of Handelman
mentioned above.

Finally, Problem~\ref{Pb:aleph1Ring} was solved
positively by the second author in~\cite{Wehr3}.

\medskip
In \cite[Theorem 1.5]{GoHa86}, Handelman and the first author
showed that for every dimension group
$G$ of size at most $\aleph_1$, there exists a locally
matricial algebra $R$ such that $V(R)\cong G^+$ (in fact,
the result is given there in the case where
$G$ has an order-unit. In the general case, $G$ embeds as an
ideal into a dimension group
$H$ with order-unit such that $|H|\leq\aleph_1$---take, for
example, $H=\QQ\times_{\mathrm{lex}}G$, the lexicographical
product of $\QQ$ by $G$---and then we can use
Proposition~\ref{P:IdVRIR}.(i)). Therefore, by
Proposition~\ref{P:ManyIso}, the analogue of
Problem~\ref{Pb:aleph1Ring} for locally matricial algebras
(\emph{i.e.}, the question whether the $\aleph_1$ version of
Bergman's Theorem holds) is equivalent to
Problem~\ref{Pb:GrpSemil}$'$.

\begin{problem}\label{Pb:DistrRing}
Let $S$ be a distributive $0$-\emph{lattice}. Does there exist
a regular ring $R$ such that $S\cong\Idc R$?
\end{problem}

By Theorem~\ref{T:LiftDLat}, every distributive $0$-lattice
is isomorphic to the maximal semilattice quotient of $G^+$
for some dimension group $G$. However, this does not help
because there are dimension groups of size
$\aleph_2$ that are not isomorphic to $K_0(R)$ for any regular
ring $R$ (see \cite{Wehr1}).

Finally, P. R\r u\v zi\v cka solved Problem~\ref{Pb:DistrRing}
positively in \cite{Ruzi1}.

\medskip

Natural extensions of the problems above are found when one
does not just ask for lifting  semilattices, but their
\emph{homomorphisms}. The solution of lattice-theoretical
analogues of this kind of problem can be found in
\cite{Wehr4,Wehr5}.

\begin{problem}\label{Pb:Id1}
Characterize the distributive $0$-semilattices $S$ such that
for every locally matricial algebra $R$, every
$\{\vee,0\}$-homomorphism $\varphi\colon\Idc R\to S$
can be lifted, that is, there are a locally matricial algebra
$R'$, an algebra homomorphism $f\colon R\to R'$, and an isomorphism
$\alpha\colon\Idc R'\to S$ such that
$\alpha\circ\Idc f=\varphi$.
\end{problem}

Of course, the map $\Idc f$ is defined by the rule
 \[
 (\Idc f)(xR)=f(x)R',\qquad\text{for all }x\in R,
 \]
thus turning $\Idc$ into a \emph{functor}.
\medskip

Pursuing the lattice-theoretical analogy, it is reasonable to
ask for the following two-dimensional analogue of
Problem~\ref{Pb:Id1}:

\begin{problem}\label{Pb:Id2}
Characterize the distributive $0$-semilattices $S$ such that
for every diagram $\mathcal{D}$ of locally matricial algebras
of the form $f_i\colon R\to R_i$, for $i\in\{1,2\}$, every
homomorphism $\varphi\colon\Idc\mathcal{D}\to S$ can be
lifted by some homomorphism $f\colon\mathcal{D}\to R'$, for
some locally matricial algebra $R'$.
\end{problem}

The second author's paper \cite{Wehr5} studies the
$0$-semilattices $S$ that satisfy a
lattice-theoretical analogue of Problem~\ref{Pb:Id2}.


\begin{thebibliography}{99}

\bibitem{AGOMP}
P. Ara, K.R. Goodearl, K.C. O'Meara, and E. Pardo,
\emph{Separative cancellation for projective modules over
exchange rings},
Israel J. Math. \textbf{105} (1998), 105--137.

\bibitem{Berg86}
G.M. Bergman,
\emph{Von Neumann regular rings with tailor-made ideal
lattices}, Unpublished note (26 October 1986).

\bibitem{Birk67}
G. Birkhoff,
``Lattice theory. Corr. repr. of the 1967 3rd ed.'',
Third Ed., American Math. Soc. Colloq. Publ. \textbf{25},
Providence, Rhode Island: American Mathematical Society,
vi + 418~pp. (1979).

\bibitem{Blac86}
B. Blackadar,
``K-Theory for Operator Algebras'', MSRI Publ. \textbf{5},
Springer-Verlag, New York, Springer-Verlag, 1986,
vii + 338~pp.

\bibitem{ClPr}
A.H. Clifford and G.B. Preston,
``The Algebraic Theory of Semigroups'',
Mathematical Surveys \textbf{7}, American Math. Soc.,
Providence, R.I., Vol.~\textbf{1}, 1961, xv + 224~pp., and
Vol.~\textbf{2}, 1967, xvi + 350~pp.

\bibitem{Cntz81}
J. Cuntz,
\emph{K-theory for certain C*-algebras},
Ann. of Math. \textbf{113} (1981), 181--197.

\bibitem{Dav96}
K.R. Davidson,
``C*-Algebras by Example'', Fields Inst. Monographs
\textbf{6}, American Math. Soc., Providence, 1996, xiv +
309~pp.

\bibitem{Dobb82}
H. Dobbertin,
\emph{On Vaught's criterion for
isomorphisms of countable Boolean algebras},
Algebra Universalis \textbf{15} (1982), 95--114.

\bibitem{EHS80}
E.G. Effros, D.E. Handelman and C.-L. Shen,
\emph{Dimension groups and their affine representations},
Amer. J. Math. \textbf{102} (1980), no.~2,
385--407.

\bibitem{Elli76}
G.A. Elliott,
\emph{On the classification of inductive limits of sequences of
semisimple finite-dimensional algebras},
J. Algebra \textbf{38} (1976), 29--44.

\bibitem{FrHa56}
K.D. Fryer and I. Halperin,
\emph{The von~Neumann coordinatization theorem for complemented
modular lattices},
Acta Sci. Math. (Szeged) \textbf{17} (1956),
203--249.

\bibitem{Gvnrr}
K.R. Goodearl,
``Von Neumann Regular Rings'', Pitman, London 1979, xvii +
369~pp.; Second Ed. Krieger, Malabar, Fl., 1991, xvi + 412~pp.

\bibitem{Gartnoeth}
\bysame,
\emph{Artinian and n\oe{}therian modules over regular rings},
Comm. Algebra \textbf{8} (1980), 477--504.

\bibitem{Gpoag}
\bysame,
``Partially Ordered Abelian Groups with Interpolation'',
Math. Surveys and Monographs \textbf{20}, Amer. Math. Soc.,
Providence, 1986, xxii + 336~pp.

\bibitem{GoHa86}
K.R. Goodearl and D.E. Handelman,
\emph{Tensor products of dimension groups and $K_0$ of
unit-regular rings}, Canad. J. Math. \textbf{38},
no.~3 (1986), 633--658.

\bibitem{Grat71}
G. Gr\"atzer,
``Lattice Theory. First Concepts and Distributive Lattices'',
W.~H. Freeman and Co., San Francisco, Calif., 1971.
xv + 212~pp.

\bibitem{Grat98}
\bysame,
``General Lattice Theory. Second Edition'',
Birkh\"auser Verlag, Basel, 1998, xix + 663~pp.

\bibitem{GrSc}
G. Gr\"atzer and E.T. Schmidt,
\emph{Congruence lattices of lattices},
Appendix C in \cite{Grat98}.

\bibitem{Gril76}
P.A. Grillet,
\emph{Directed colimits of free commutative semigroups},
J. Pure Appl. Algebra \textbf{9} (1976), no.~1, 73--87.

\bibitem{Hndlmn}
D.E. Handelman,
\emph{Notes on ideal lattices},
Unpublished note (1981).

\bibitem{Jac64}
N. Jacobson,
``Structure of Rings'',
Rev. Ed., Colloq. Publ. \textbf{37}, Amer. Math. Soc.,
Providence, 1964, vii + 263~pp.

\bibitem{KR}
K.H. Kim and F.W. Roush,
\emph{Regular rings and distributive lattices},
Comm. Algebra \textbf{8} (1980), 1283--1290.

\bibitem{Lin97}
H. Lin,
\emph{C*-algebras of trivial K-theory and semilattices},
Internat. J. Math. \textbf{10} (1999), 93--128.

\bibitem{Malc}
A.I. Mal'cev,
``Algebraic Systems'',
Die Grundlagen der mathematischen
Wissenschaften in Einzeldarstellungen, Band 192,
Springer-Verlag, Berlin Heidelberg New York, 1973.
xii + 317~pp.

\bibitem{Mur90}
G.J. Murphy,
``C*-Algebras and Operator Theory'',
Academic Press, Boston, 1990, x + 286~pp.

\bibitem{PaSa79}
W. Paschke and N. Salinas,
\emph{Matrix algebras over $\mathcal{O}_n$},
Michigan Math. J. \textbf{26} (1979), 3--12.

\bibitem{Pu}
P. Pudl\'ak,
\emph{On congruence lattices of lattices},
Algebra Universalis \textbf{20} (1985), 96--114.

\bibitem{Ruzi1}
P. R\r u\v zi\v cka,
\emph{Lattices of two-sided ideals of locally matricial
algebras and the $\Gamma$-invariant problem}, preprint.

\bibitem{Ruzi2}
\bysame,
\emph{A distributive semilattice that is not isomorphic to the
compact ideal lattice of any dimension group}, preprint.

\bibitem{Schm84}
E.T. Schmidt,
\emph{Congruence lattices of complemented modular lattices},
Algebra Universalis \textbf{18} (1984), 386--395.

\bibitem{Shan74}
R.T. Shannon,
\emph{Lazard's theorem in algebraic categories},
Algebra Universalis \textbf{4} (1974), 226--228.

\bibitem{ShWe94}
R.M. Shortt and F. Wehrung,
\emph{Common extensions of semigroup-valued charges},
J. Math. Anal. Appl. \textbf{187}, no.~1 (October 1,
1994), 235--258.

\bibitem{vonN36}
J. von~Neumann,
\emph{On regular rings},
Proc. Nat. Acad. Sci. USA \textbf{22} (1936), 707--713.

\bibitem{Wegg93}
N.E. Wegge-Olsen,
``K-Theory and C*-Algebras, A Friendly Approach'',
Oxford Univ. Press, Oxford, 1993, xii + 370~pp.

\bibitem{Wehr92}
F. Wehrung,
\emph{Injective positively ordered monoids I},
J. Pure Appl. Algebra \textbf{83} (1992), 43--82.

\bibitem{Wehr1}
\bysame,
\emph{Non-measurability properties of interpolation
vector spaces},
Israel J. Math. \textbf{103} (1998),
177--206.

\bibitem{Wehr2}
\bysame,
\emph{A uniform refinement property for
congruence lat\-tices},
Proc. Amer. Math. Soc. \textbf{127}, no.~2
(1999), 363--370.

\bibitem{Wehr3}
\bysame,
\emph{Representation of algebraic distributive lattices
with $\aleph_1$ compact elements as ideal lattices of regular
rings}, Publ. Mat., to appear.

\bibitem{Wehr5}
\bysame,
\emph{Forcing extensions of partial lattices}, preprint.

\bibitem{Wehr4}
\bysame,
\emph{Semilattices with the two-dimensional congruence
amalgamation property}, preprint.

\bibitem{Zha90}
S. Zhang,
\emph{A Riesz decomposition property
and ideal structure of multiplier algebras},
J. Operator Theory \textbf{24} (1990), 204--225.

\end{thebibliography}
\end{document}